\documentclass[numbook, envcountsect, envcountsame, envcountreset,runningheads, smallextended]{svjour3}

\smartqed  

\usepackage{graphicx}
\usepackage{fix-cm}

\usepackage{amsmath}
\usepackage{amsfonts}
\usepackage{amssymb}

\def\JELname{\textbf{Journal of Economic Literature Classification}\enspace}
\def\JEL#1{\par\addvspace\medskipamount{\rightskip=0pt plus1cm
\def\and{\ifhmode\unskip\nobreak\fi\ $\cdot$
}\noindent\JELname\ignorespaces#1\par}}

\newtheorem{condition}[theorem]{Condition}

\newcommand{\re}{{\mathbb R}}    

\newcommand{\Var}{\operatorname{Var}}

\newcommand{\PolmS}{\mathcal{P}_m}

\newcommand{\im}{\operatorname{i}}

\begin{document}

\title{Polynomial processes and their applications to mathematical finance}

\titlerunning{Polynomial processes and applications to mathematical finance}

\author{Christa Cuchiero \and Martin Keller-Ressel \and Josef Teichmann}
\authorrunning{C.~Cuchiero, M.~Keller-Ressel, J.~Teichmann}

\institute{Christa Cuchiero \at
           University of Vienna, Faculty of Mathematics, Nordbergstrasse 15, A-1090 Wien, Austria \\
           Tel.:+43 1 4277 50757\\
           \email{christa.cuchiero@univie.ac.at}
           \and
           Martin Keller-Ressel \at
           TU Berlin - Fakult\"at II, Institut f\"ur Mathematik, Strasse des 17. Juni 136, D-10623 Berlin, Germany \\
           Tel.:+49 30 314 28619\\
           \email{mkeller@math.TU-Berlin.DE}
           \and
           Josef Teichmann \at
           ETH Z\"urich, D-MATH, R\"amistrasse 101, CH-8092 Z\"urich, Switzerland \\
           Tel.:+41 44 632 3174\\
           \email{josef.teichmann@math.ethz.ch}
}

\date{Received: date / Accepted: date}

\maketitle

\begin{abstract}
We introduce a class of Markov processes, called $m$-polynomial, for which the calculation of (mixed) moments up to order $m$ only requires the computation of matrix exponentials. This class contains affine processes, processes with quadratic diffusion coefficients, as well as L\'evy-driven SDEs with affine vector fields. Thus, many popular models such as exponential L\'evy models or affine models are covered by this setting. The applications range from statistical GMM estimation procedures to new techniques for option pricing and hedging. For instance, the efficient and easy computation of moments can be used for variance reduction techniques in Monte Carlo methods.

 \JEL{C02, G12}
\end{abstract}

\section{Introduction}\label{intro}

Pricing and hedging of contingent claims are the crucial computations done within every model in mathematical finance. For European type claims this amounts to the computation of the expected value of a functional of the (discounted) price process under some martingale measure. (Partial) Hedging portfolios are then constructed via appropriate derivatives of those expected values with respect to model parameters or to the prices, so called Greeks. Let us denote the (discounted) price process at time $T$, a vector in $ \mathbb{R}^n $, by $ X_T $. We can roughly distinguish three cases of complexity for the mentioned computations:

\begin{enumerate}

\item The probability distribution of $ X_T $ is analytically known.
\item The characteristic function of $ X_T $ is analytically known.
\item The local characteristics of $ X_T $ are analytically known.
\end{enumerate}

In the first case, a numerical quadrature algorithm is sufficient for the efficient computation of the contingent claim's price $ \mathbb{E}[\phi(X_T)] $, where $\phi$ denotes some payoff function.

In the second case, variants of Plancherel's theorem are applied in order to evaluate the price functional $\mathbb{E}[\phi(X_T)] $, for instance,
\[
\mathbb{E}[\phi(X_T)] = \int_{\mathbb{R}^n} \widehat{\phi}(u) \mathbb{E}[ \exp( i\langle u , X_T \rangle)] du,
\]
where $ \widehat{\phi} $ denotes the Fourier transform of the function $ \phi $. Remark that often modifications of the original payoff function are used to make the Fourier methodology applicable. This is numerically efficient, even though its implementation, in particular the complex integration, can take some time (see, e.g.,~\cite{carr}). Also there are different levels of what it means to ``know analytically'' the characteristic function of $ X_T $. In affine models, for example, one might need to solve a high-dimensional Riccati equation for each $ u \in \mathbb{R}^n $ to calculate the characteristic function $ u \mapsto \mathbb{E}[ \exp( i\langle u , X_T \rangle)] $. In this case, ``analytic knowledge'' involves some precalculations, which also have to be performed efficiently. 

The third case is characterized by the use of Monte Carlo methods: one samples from the (unknown) distribution of $ X_T $ by generating, for instance through Euler schemes, approximate distributions for $ X_T $. This procedure is very robust, but takes a considerable amount of time. Moreover, for the convergence of the Euler scheme certain regularity assumptions on the characteristics, e.g.~(local) Lipschitz continuity, are required.

In this article we would like to add a fourth case which -- in the previous order -- would correspond to case $ 1 \frac{1}{2} $. 
We can describe a class of Markov processes, called ``polynomial processes'', which have the property that the expected 
value of \emph{any polynomial} of the random variables $X_t$, $t \geq 0$, is \emph{again a polynomial in the initial value} of the process. This means in particular that
moments of all orders of $ X_T $ can be computed in an \emph{easy and efficient} way, even though neither its probability distribution nor its characteristic function needs to be known. Loosely speaking one could say that the expressions for all finite moments are analytically known (up to a matrix exponential). 

We shall analyze this class and show that exponential L\'evy processes, affine processes or processes of Pearson diffusion type belong to it. 
The method is best explained by an example: consider a stochastic volatility model of SVJJ-type~\cite{gatheral}, 
i.e.~both the logarithmic (discounted) price process and the stochastic volatility can jump. Such models can be described by stochastic differential equations of the type
\begin{align*}
dY_t & = (b_1+\beta_{11}Y_t+\beta_{21}V_t)dt + \sqrt{V_t} dB_{t,1} + dZ_{t,1}, \\
dV_t& = (b_2+\beta_{22}V_t) dt + \sqrt{V_t} dB_{t,2} + dZ_{t,2},
\end{align*}
where $ (B_1,B_2) $ are possibly correlated Brownian motions and $ (Z_1,Z_2) $ is a bivariate pure-jump L\'evy process, independent of $ (B_1,B_2) $, 
whose second component $ Z_2 $ has positive increments. For such models there is no easy-to-implement (explicit) formula for the characteristic function, 
even though they are affine models. Assuming now appropriate moment conditions on the jump measures, the Markov process $ X=(Y,V) $ turns out to be a 
polynomial process, that is, the expected value of any polynomial of $(Y_t,V_t) $ is a polynomial in $ Y_0 $ and $ V_0 $. 
The coefficients of this polynomial can be calculated efficiently by exponentiating a matrix, which can be easily deduced from the (extended) generator. In other words, there is a large subset of claims for which the prices and hedge ratios are explicitly known (up to matrix exponentials). Large can be made precise in the following sense: if the law of $ X_T $, say $\mu$, is characterized by its moments, then ``large'' means dense, i.e.~polynomial claims are dense (with respect to the $L^1(\mu)$ norm) in the set of ``all''  claims. If this is not the case, the payoff function can at least be uniformly approximated by polynomials on some interval, which can be chosen according to the support of the probability distribution. The explicit knowledge of prices of polynomial claims then allows to apply variance reduction techniques for Monte-Carlo computations. 

The remainder of our article is organized as follows: in Section~\ref{sec-pol-proc} we formally introduce the class of $m$-polynomial processes, establish a relationship to semimartingales and give conditions on the (extended) generator such that a Markov process is $m$-polynomial. Section~\ref{sec:examples} deals with examples from the class of $m$-polynomial processes and Section~\ref{sec-applications} with applications to pricing and hedging in mathematical finance.

\section{Polynomial Processes}\label{sec-pol-proc}

We define polynomial processes as a particular class of time-homogeneous Markov processes with state space $S \subseteq \mathbb{R}^n$, 
some closed subset of $\re^n$. To clarify notation, we find it useful to recall the basic ingredients of a time-homogeneous Markov process 
and the particular assumptions and conventions being made in this article (compare~\cite[Chapter 3]{rogers}).

Throughout, $S$ is a closed subset of $\re^n$ and $\mathcal{S}$ denotes its Borel $\sigma$-algebra.
Since we shall not assume the process to be conservative, we adjoin to the state space $S$ a point $\Delta \notin S$, called cemetery, 
and set $S_{\Delta}=S\cup\{\Delta\}$ as well as $\mathcal{S}_{\Delta}=\sigma(\mathcal{S}, \{\Delta\})$.
We make the convention that $f(\Delta)=0$ for any function $f$ on $S$.

We consider a time-homogeneous Markov semigroup $(P_t)_{t \geq 0}$ given by
\[
 P_tf(x):=\int_{S} f(\xi)p_t(x,d\xi), \quad x \in S_{\Delta},
\]
and acting on all Borel measurable functions $f: S_{\Delta} \rightarrow \re$ for which the integral is well defined. 
Here, $(p_t)_{t \geq 0}$ denotes the transition function, which satisfies beside the standard conditions (see~\cite[Definition III.1.1]{rogers})
the following properties:
\begin{itemize}
 \item[(i)]for all $x \in S_{\Delta}$, $p_0(x,\cdot)=\delta_{x}$, where $\delta_{x}$ denotes the Dirac measure;
\item[(ii)] for all $t\geq 0$ and $x\in S$, $p_t(x,\{\Delta\})=1-p_t(x,S)$ and $p_t(\Delta,\{\Delta\})=1$.
\end{itemize}

Since the theory of polynomial processes always deals with a Markov process $(X_t)_{t\geq 0}$ having the property that $(f(X_t))_{t\geq 0}$ is a special semimartingale for all linear functionals $f$ on $\mathbb{R}^n$ (extended by $0$ on $\Delta$), we can bona fide assume that the probability space $\Omega$ is the space of c\`adl\`ag functions $\omega: \mathbb{R}_+ \to S_{\Delta}$ such that $\omega(t-)=\Delta$ and $\omega(t)=\Delta$ implies $\omega(s)=\Delta$ for all $s \geq t$.  This follows from the simple conclusion (i) to (iii) of Theorem \ref{th:characterizationpoly}, martingale regularity, and from the remarks at the bottom of \cite[p.~245]{rogers}.
We thus understand $X$ as the coordinate process $X_t(\omega)=\omega(t)$ and denote by $(\mathcal{F}_t^0)$ the filtration generated by $X$ and set $\mathcal{F}^0=\bigvee_{t \geq 0}\mathcal{F}_t^0$.
 
In the sequel we consider some right continuous filtration $(\mathcal{F}_t)_{t \geq 0}$ satisfying $\mathcal{F}_t^0 \subseteq \mathcal{F}_t$ and 
$\mathcal{F}=\bigvee_{t \geq 0}\mathcal{F}_t$. 
We finally assume that for each $x \in S_{\Delta}$ there exists a probability measure $\mathbb{P}_x$ on $(\Omega, \mathcal{F})$ such that $X$ is Markovian relative 
to $(\mathcal{F}_t)$ with semigroup $(P_t)$, that is,
\[
\mathbb{E}_x[f(X_{t+s})|\mathcal{F}_s]=\mathbb{E}_{X_s}[f(X_{t})]=P_tf(X_s), \quad \mathbb{P}_x\textrm{-a.s.}
\]
for all $x \in S_{\Delta}$, $s,t \in [0,\infty)$ and all Borel functions $f: S_{\Delta} \rightarrow \re$ satisfying\\ 
$\mathbb{E}_x[|f(X_t)|]< \infty$ for all $t \geq 0$ and $x \in S$. 

\subsection{Definition and characterization of polynomial processes}\label{sec:2.1}

For the treatment of polynomial process we have chosen a framework of stochastic analysis leading to easy-to-verify conditions for a process to be $m$-polynomial (see Theorem~\ref{th:characterizationpoly(iv)} and Theorem~\ref{cor:suffcond}). 
Before giving the precise definition of polynomial processes, let us introduce some notation.

Let $\PolmS$ denote the finite dimensional vector space of polynomials up to degree $m \geq 0$ on $S$, i.e.~the restriction of polynomials on $\re^n$ to $S$, defined by
\begin{align*}
\PolmS:=\left\{S \ni x \mapsto \sum_{|\mathbf{k}|=0}^m \alpha_{\mathbf{k}}x^{\mathbf{k}},\, \Delta \mapsto 0 \,\Big |\, \alpha_{\mathbf{k}} \in \re\right\},
\end{align*}
where we use multi-index notation $\mathbf{k} = (k_1,\ldots, k_n) \in \mathbb{N}^n_0$, $|\mathbf{k}|=k_1+\cdots+k_n$ and $x^{\mathbf{k}}=x^{k_1}_1\cdots x^{k_n}_n$. The dimension of $\PolmS$ is denoted by $N < \infty$ and depends on $ S $: if $ S $ is a single point, the dimension is always $ 1 $ and if $ S $ is the whole space $ \mathbb{R}^n $, it is maximal.

Moreover, for every multi-index $\mathbf{k}$, we define functions $f_{\mathbf{k}}$ by setting
\begin{align}\label{not:f}
f_{\mathbf{k}}=\left \{\begin{array}{ll}x^\mathbf{k}&\textrm{if } x \in S,\\
0&\textrm{if } x=\Delta.
\end{array} \right.
\end{align}
Furthermore, we write $f_i$ when $\mathbf{k}=e_i$ and $f_{ij}$ when $\mathbf{k}=e_i+e_j$, where $e_i,\, i\in \{1, \ldots, n\},$ denotes the $i^{\textrm{th}}$ canonical basis vector.
Then the space $\mathcal{P}_m$ clearly corresponds to the linear hull of the functions $\{f_{\mathbf{k}}, \, |\mathbf{k}|\leq m\}$.

Here is our main definition:

\begin{definition}\label{def-m-pol}
We call an $S_{\Delta}$-valued time-homogeneous Markov process\\
\emph{$m$-polynomial} if we have for all $k\in\{0,\ldots, m\}$, all $f \in \mathcal{P}_k$, $x \in S$ and $t \geq 0$,
\[
x \mapsto P_tf(x) \in \mathcal{P}_k.
\]
Additionally, we assume that $t \mapsto P_tf(x)$ is continuous at $t=0$ for all $f \in \PolmS$.
If $X$ is $m$-polynomial for all $m \geq 0$, then it is called \emph{polynomial}.
\end{definition}

\begin{remark}
\begin{itemize}
\item[(i)] Let us stress that in the above definition it is implicitly assumed that
\[
P_t|f|(x)=\mathbb{E}_x[|f(X_t)|]=\int_{S} |f(\xi)| p_t(x,d\xi)< \infty
\]
for every  $f \in \PolmS$, $x \in S_{\Delta}$ and $t \geq 0$, because otherwise the expression 
$P_tf(x)=\mathbb{E}_x[f(X_t)]$ would not even be well-defined.
\item[(ii)] The subtlety of Definition \ref{def-m-pol} lies in the fact that we assume $ P_t (\mathcal{P}_k)\subset \mathcal{P}_k$ for \emph{all} $k \in \{0, \ldots, m\} $ (compare with Remark~\ref{rem:remarkspol} (iv)). The assumption that $ P_t (\PolmS) \subset \PolmS $ only for $m$, but not for smaller degrees is not sufficient for our proofs of Theorem~\ref{th:characterizationpoly(iv)} and Theorem~\ref{cor:suffcond}, which we consider as the most important assertions from the point of view of applications.
\end{itemize}
\end{remark}

Let us introduce the following notion of an extended Markov generator, which is due to Dynkin (see, e.g.,~\cite[Definition 7.1]{cinlar}) 
and which we shall use to characterize $m$-polynomial processes.

\begin{definition}\label{def:extendedgenerator}
An operator $\mathcal{G}$ with domain $\mathcal{D}_{\mathcal{G}}$ is called \emph{extended generator} for some Markov process $X$, if $\mathcal{D}_{\mathcal G}$ consists of those Borel measurable functions $f: S_{\Delta} \rightarrow \mathbb{R}$ for which there exists a function $\mathcal{G}f$ such that the process
\begin{align}\label{eq:extendedgenerator}
M_t^f:=f(X_{t })-f(x)-\int_0^{t } \mathcal{G}f(X_s) ds
\end{align}
is well defined and a $(\mathcal{F}_t,\mathbb{P}_x)$-local martingale for every $x \in S_{\Delta}$.
\end{definition}

\begin{remark}\label{rm:localmartingale}
We define the lifetime of the process by
\begin{align}\label{eq:lifetime}
 T_{\Delta}(\omega)=\inf\{t\,|\, X_t(\omega)=\Delta\},
\end{align}
where the infimum over the empty set is set to be $\infty$. Since
$
 \{ T_{\Delta} < t \}= \bigcup_{q < t, q \in \mathbb{Q}} \{X_q =\Delta\} \in \mathcal{F}_t
$
and as $(\mathcal{F}_t)$ is supposed to be right continuous, $T_{\Delta}$ is an $\mathcal{F}_t$-stopping time.
Due to our convention $f(\Delta)=0$, the local martingale property of~\eqref{eq:extendedgenerator} is therefore equivalent to
\[
 f(X_t)1_{\{t < T_{\Delta}\}}-f(x)- \int_0^{t \wedge T_{\Delta}} \mathcal{G}f(X_s) ds
\]
being a local martingale. 
\end{remark}

\begin{remark}\label{remarks_def:generator}
Suppose that $f$ lies in the domain of the extended generator and satisfies $P_t|f|(x)< \infty$ for all $t \geq 0$ and $x \in S$.
Then $M^f$ as defined in~\eqref{eq:extendedgenerator} is a \emph{true} martingale if and only if all increments of $ f(X_t)-f(x)- \int_0^t \mathcal{G}f(X_s) ds $ have vanishing expectation, i.e.~for all $u < t$,
\[
\mathbb{E}_x\left[f(X_t)-f(X_u)-\int_u^t \mathcal{G}f(X_s)ds\right]=P_{t}f(x)-P_uf(x)-\int_u^t P_s\mathcal{G}f(x)ds=0.
\]
In particular, by Fubini's theorem $ \int_0^t P_s \mathcal{G} f(x) ds $ exists on finite time intervals and thus also
$ P_s | \mathcal{G} f |(x) $ for almost all $ s $ with respect to the Lebesgue measure. 
\end{remark}

The following lemma is well known for Feller processes (but perhaps not in our particular setting) and aims to
establish a connection between the Kolmogorov backward equation and the extended generator introduced in~\eqref{def:extendedgenerator}.

\begin{lemma}\label{lem:Ptf}
Let $X$ be a time-homogeneous Markov process with semigroup $(P_t)$ and denote by $f:S_{\Delta} \to \mathbb{R}$ some function satisfying $P_t|f|(x)< \infty$ for all $t \geq 0$ and $x \in S$.
If $f$ lies in the domain of the extended generator, $f \in \mathcal{D}_{\mathcal{G}}$, and if $M^f$ as defined in~\eqref{eq:extendedgenerator} 
is additionally a true martingale, then we have:
\begin{itemize}
\item[(i)] For any given $s \geq 0$, $M^{P_sf}$ is a true martingale, in particular 
$P_sf \in \mathcal{D}_{\mathcal{G}}$, and $\mathcal{G}P_sf=P_s\mathcal{G}f$.
\item[(ii)] If $t \mapsto P_t \mathcal{G}f(x)$ is continuous at $t=0$, then $P_tf$ solves the Kolmogorov backward equation, that is,
\[
\frac{\partial u (t, x)}{\partial t}=\mathcal{G}u(t, x), \quad u(0, x)=f(x).
\]
\end{itemize}
\end{lemma}

\begin{proof}
For the first statement, we show that
\[
P_sf(X_t)-P_sf(x)- \int_0^t P_s\mathcal{G}f(X_r) dr
\]
is a true $(\mathcal{F}_t,\mathbb{P}_x)$-martingale for any fixed $s \geq 0$. 
By the definition of the extended generator, this then implies that $P_sf \in \mathcal{D}_{\mathcal{G}}$ and $\mathcal{G}P_sf=P_s\mathcal{G}f$.
Indeed, we have by the assumption $P_t|f|(x)< \infty$ for all $t \geq 0$ and $x \in S$ and Remark~\ref{remarks_def:generator} that 
$f(X_t)$ and $\mathcal{G}f(X_t)$ are integrable for every $t \geq 0$, hence $P_sf(X_t)$ and $P_s\mathcal{G}f(X_t)$ as well. 
Therefore the following expectation is well defined and we obtain for $u \leq t$
\begin{align*}
&\mathbb{E}_x\left[P_sf(X_t)-P_sf(x)- \int_0^t P_s\mathcal{G}f(X_r) dr\,\Big|\,\mathcal{F}_u\right]\\
&\quad=P_sf(X_u)-P_sf(x)- \int_0^u P_s\mathcal{G}f(X_r) dr\\
&\quad \quad + \mathbb{E}_x\left[P_sf(X_t)-P_sf(X_u)- \int_u^t P_s\mathcal{G}f(X_r) dr\,\Big|\,\mathcal{F}_u\right].
\end{align*}
By the Markov property, the conditional expectation on the right is equal to
\[
\mathbb{E}_{X_u}\left[P_sf(X_{t-u})-P_sf(X_0)- \int_0^{t-u} P_s\mathcal{G}f(X_r)dr\right].
\]
But for any $y \in S_{\Delta}$, we have
\begin{align*}
&\mathbb{E}_{y}\left[P_sf(X_{t-u})-P_sf(X_0)- \int_0^{t-u} P_s\mathcal{G}f(X_r)dr\right]\\
&\quad =P_{s+t-u}f(y)-P_sf(y)-\int_s^{s+t-u} P_{\tilde{r}}\mathcal{G}f(y)d \tilde{r}\\
&\quad =0,
\end{align*}
where the last equality follows from Remark~\ref{remarks_def:generator}. This completes the proof of (i).

Statement (ii) follows from  Remark~\ref{remarks_def:generator}, the continuity of $t \mapsto P_t\mathcal{G}f(x)$ and from assertion (i), since
\begin{align*}
\frac{\partial P_tf(x)}{\partial t}&= \lim_{h \rightarrow 0} \frac{P_{t+h}f(x)-P_tf(x)}{h}=\lim_{h \rightarrow 0} P_t\frac{P_{h}f(x)-f(x)}{h}\\
&=\lim_{h \rightarrow 0} P_t\frac{1}{h}\int_0^h P_s\mathcal{G}f(x)ds=P_t\mathcal{G}f(x)=\mathcal{G}P_tf(x).
\end{align*}
\end{proof}

Let us now state our first theorem which is a consequence of elementary results in semigroup theory:

\begin{theorem}\label{th:characterizationpoly}
Let $X$ be a time-homogeneous Markov process with state space $S_{\Delta}$ and semigroup $(P_t)$. 
Then the following assertions are equivalent:
\begin{itemize}
\item[(i)] $X$ is $m$-polynomial for some $m \geq 0$.
\item[(ii)] For every $k \in \{0, \ldots, m\}$, there exists a linear map $A$ on $ \mathcal{P}_k$, such that for all $t \geq 0$, 
$(P_t)$ restricted to $ \mathcal{P}_k$ can be written as
\[
P_t|_{ \mathcal{P}_k}=e^{tA}.
\]
\item[(iii)] 
For all $f \in \PolmS$, $x \in S_{\Delta}$ and $t \geq 0$, $P_t|f|(x)=\mathbb{E}_x[|f(X_t)|]< \infty$ and
\[
M_t^f=f(X_{t})-f(x)-\int_0^{t} \mathcal{G} f (X_s) ds
\]
is a $(\mathcal{F}_t,\mathbb{P}_x)$ true martingale, and the extended generator $\mathcal{G}$ satisfies\\ 
$\mathcal{G}(\mathcal{P}_k) \subset \mathcal{P}_k$ for all $k \in \{0,1\ldots, m\}$.
\end{itemize}
\end{theorem}

\begin{proof}
Our strategy to prove the above equivalences is to show (i) $\Rightarrow$ (ii) $\Rightarrow$ (iii) $\Rightarrow$ (i).

Throughout the proof, let $k \in \{0, \ldots, m\}$ be fixed.
We start by showing (i) $\Rightarrow$ (ii). 
By the definition of an $m$-polynomial process, the Markovian semigroup $(P_t)$ induces a semigroup of operators on $\mathcal{P}_k$.
Since $ \mathcal{P}_k$ is finite dimensional and $t \mapsto P_{t}|_{ \mathcal{P}_k}$ continuous at $t=0$, standard results of semigroup theory (see, e.g.,~\cite[Theorem 2.9]{nagel}), imply the representation of $P_{t}|_{ \mathcal{P}_k}$ as matrix exponential, that is, there exists some linear map $A$ on $\mathcal{P}_k$ such that $P_t|_{ \mathcal{P}_k}=e^{tA}$.

Next, we show (ii) $\Rightarrow$ (iii). For every $f \in \mathcal{P}_k$ we have by (ii), $Af \in \mathcal{P}_k$ and
\[
P_tf-f-\int_0^tP_sAf ds=e^{tA}f-f-\int_0^t e^{sA}A f ds=0.
\]
Thus,
$
f(X_t)-f(x)- \int_0^t Af(X_s) ds
$
is a $(\mathcal{F}_t,\mathbb{P}_x)$-martingale. Hence $f$ lies in $D_{\mathcal{G}}$ and $\mathcal{G}f = Af$, implying that $\mathcal{G}(\mathcal{P}_k) \subset \mathcal{P}_k$ holds true.

In order to prove (iii) $\Rightarrow$ (i), we consider the Kolmogorov backward equation for an initial value $u(0,\cdot)=f \in \mathcal{P}_k$
\[
\frac{\partial u(t, x)}{\partial t}=\mathcal{G}u(t,x).
\]
By Lemma~\ref{lem:Ptf} (ii), $P_tf$ solves the Kolmogorov equation, since $t \mapsto P_t\mathcal{G}f(x)$ is continuous at $t=0$ for any $f \in \mathcal{P}_k$. 
This follows from the fact that $\mathcal{G}$ maps $\mathcal{P}_k$ to itself and the martingale property of $M^f$, which implies
\[
P_tf(x)=f(x)+\int_0^t P_s\mathcal{G}f(x) ds,
\]
thus in particular continuity of $t \mapsto P_tf(x)$ for any $f \in \mathcal{P}_k$. 
By choosing a basis $\langle e_1, \ldots, e_N \rangle$ of $\mathcal{P}_k$, we can define a linear map $A$ on  $\mathcal{P}_k$ by setting
\[
\mathcal{G}e_i =:\sum_{j=1}^N A_{ij} e_j.
\]
Then $\mathcal{G}|_{\mathcal{P}_k}=A$. On $\mathcal{P}_k$, the Kolmogorov backward equation thus
reduces to the following linear ODE
\[
\frac{\partial u(t)}{\partial t}=Au(t), \quad u(0)=f,
\]
whose unique solution is given by $e^{tA}f$ (see, e.g., ~\cite[Theorem 2.8]{nagel}).
Hence on $\mathcal{P}_k$, $P_tf$  is equal to $e^{tA}f$ and is therefore a polynomial of degree smaller than or equal to $k$. 
Since this holds true for all $k \in \{0,1, \ldots,m\}$, $X$ is $m$-polynomial.
\end{proof}

\begin{remark}
There is no need in assertion~(ii) to restrict the time parameter $t$ to $\re_+$, since for $t \in \re$, $(e^{tA})$ extends to a group.
\end{remark}

The equivalence (iii) $\Leftrightarrow$ (i) of Theorem~\ref{th:characterizationpoly} provides a characterization of $m$-polynomial processes in terms of the extended generator, however under the additional 
assumption that $M^f$ as defined in~\eqref{eq:extendedgenerator} is a \emph{true} martingale. If $m \geq 2$ is an even number, this latter condition is no longer needed, since 
the local martingales $M^f$ turn out to be always true martingales. 
In other words, for even numbers $m \geq 2$, Condition~\eqref{eq:cond_extended_gen} below is necessary and sufficient for $X$ being an $m$-polynomial process.

\begin{condition}\label{eq:cond_extended_gen}
$\PolmS$ lies in the domain of the extended generator, i.e.~for all $f \in \PolmS$, $x \in S_{\Delta}$ and $t \geq 0$,
$M^f$ as defined in ~\eqref{eq:extendedgenerator} is a $(\mathcal{F}_t,\mathbb{P}_x)$-\emph{local} martingale, and $\mathcal{G}(\mathcal{P}_k) \subset \mathcal{P}_k$ for all $k \in \{0,1, \ldots,m\}$
\end{condition}

\begin{theorem}\label{th:characterizationpoly(iv)}
Let $X$ be a time-homogeneous Markov process with state space $S_{\Delta}$ and let $m \geq 2$ be an even number.
Then $X$ is an $m$-polynomial process if and only if Condition~\ref{eq:cond_extended_gen} is satisfied.
\end{theorem}

\begin{proof}
The necessary direction, that is, Condition~\ref{eq:cond_extended_gen} holds, if $X$ is an $m$-polynomial process, is an obvious implication
of Theorem~\ref{th:characterizationpoly} (i) $\Rightarrow$ (iii) (of course for all $m$, not only even numbers larger than or equal to $2$).

For the sufficient direction we prove, for every $f \in \mathcal{P}_m$, $x \in S_{\Delta}$ and $t \geq 0$, that $P_t|f|(x)=\mathbb{E}_x[|f(X_t)|]< \infty$ and that $(M_t^f)$ is a true $(\mathcal{F}_t,\mathbb{P}_x)$-martingale. Then Theorem~\ref{th:characterizationpoly} (iii) $\Rightarrow$ (i) yields the assertion.

First, fix some $T > 0$ and some increasing sequence $(T_j)_{j \in \mathbb{N}}$ of stopping times with $\lim_{j \to \infty} T_j=\infty$  $\mathbb{P}_x$-a.s. such that
$ {(M^f_{t \wedge T_j})}_{t \geq 0} $ are martingales for all $f \in \mathcal{P}_m$. Furthermore we set 
\[
 F(x):=1+\sum_{i=1}^n f_i(x)^m,
\]
where $f_i, \, i \in \{1,\ldots, n\}$, are given by~\eqref{not:f}. We notice that there is some finite constant $K$ such that 
\[
 |\mathcal{G}f(x)| \leq KF(x)
\]
for all $x \in S_{\Delta}$. Hence we obtain for $0 \leq t \leq T$ and $x \in S_{\Delta}$
\begin{align*}
 \mathbb{E}_x\left[F(X_{t \wedge T_j})\right]&= F(x)+\mathbb{E}_x\left[\int_0^{t\wedge T_j}\mathcal{G}F(X_u) du\right]\\
&\leq F(x)+K\mathbb{E}_x\left[\int_0^{t\wedge T_j}F(X_u) du\right]\\
&\leq F(x)+K\mathbb{E}_x\left[\int_0^{t}F(X_{u \wedge T_j}) du\right].
\end{align*}
Since the stopping times $(T_j)$ can be chosen such that $X_{u \wedge {T_j-}}$ is bounded, 
the right-hand side of the above inequality is finite and Gronwall's lemma yields
\begin{align}\label{eq:gronwall}
 \mathbb{E}_x\left[F(X_{t \wedge T_j})\right]\leq F(x)e^{Kt}
\end{align}
for all $0 \leq t \leq T$, $j \in \mathbb{N}$ and $x \in S_{\Delta}$.
Due to the nonnegativity of $F$, we have by Fatou's lemma
\begin{align}\label{eq:finite_moment}
 \mathbb{E}_x\left[F(X_{t})\right]=\mathbb{E}_x\left[\lim_{j \to \infty}F(X_{t\wedge T_j})\right]\leq 
\liminf_{j \to \infty}\mathbb{E}_x\left[F(X_{t \wedge T_j})\right]\leq F(x)e^{Kt}, 
\end{align}
for all $0 \leq t \leq T$. Hence $P_t|f|(x)=\mathbb{E}_x[|f(X_t)|]< \infty$ for $ f \in \mathcal{P}_m $.

Next we show that for each $f \in \mathcal{P}_m$, $x \in S_{\Delta}$ and
\begin{align}\label{eq:sup}
\mathbb{E}_x\left[\sup_{t\leq T}\left|M_t^f\right|\right]< \infty.
\end{align}
First, let $f \in \mathcal{P}_k$ for $k < m$ be fixed and set $p=m/k$. Notice that there is a finite constant $\widetilde{K}$ such that
\[
 |f(x)|^p\leq \widetilde{K}F(x) \textrm{ and }  |\mathcal{G}f(x)|^p\leq \widetilde{K}F(x)
\]
for all $x \in S_{\Delta}$. We then have
\begin{align*}
\left |M^f_{t \wedge T_j}\right|^p&=\left|f(X_{t\wedge T_j})-f(x)-\int_0^{t \wedge T_j} \mathcal{G} f (X_u) du\right|^p\\
&\leq C\left(F(X_{t\wedge T_j})+F(x)+\int_0^{t \wedge T_j} F(X_u) du \right)\\
&\leq C\left(F(X_{t\wedge T_j})+F(x)+\int_0^{t} F(X_{u \wedge T_j}) du \right),
\end{align*}
for appropriate positive constants $C$ and $ t \leq T $. Taking expectations and using~\eqref{eq:gronwall}, we see that for each fixed $x$ there exists some finite constant $C_x$ such that
\begin{align*}
 \mathbb{E}_x\left[\left|M^f_{t \wedge T_j}\right|^p\right]\leq C_x,
\end{align*}
for all $j\in \mathbb{N}$ and $t \leq T$. Moreover, by Doob's maximal $L^p$-inequality for $p >1$, we have for all $j$ that
$ \mathbb{E}_x\left[\sup_{t\leq T}\left|M^f_{t \wedge T_j}\right|^p\right]\leq C \mathbb{E}_x\left[\left|M^f_{T \wedge T_j}\right|^p\right]\leq C_x$.
Since the left-hand side is increasing in $j$, monotone convergence yields~\eqref{eq:sup} for this $f$, in particular
\begin{align}\label{eq:Lp}
\sup_{t\leq T}\left|M^f_{t}\right| \in L^p.
\end{align}
Finally let us deal with the case $k=m$. We consider the polynomial $f(x)=f_i(x)^q$ for $i \in \{1, \ldots, n\}$ and $q =m/2$, which is an integer by hypothesis.
For notational convenience, we write $N=M^f$ and estimate again
\begin{align*}
f(X_t)^2&=\left(N_t+f(x)+\int_0^t \mathcal{G}f(X_u)du\right)^2\\
&\leq C\left((N_t)^2+f(x)^2+\int_0^t |\mathcal{G}f(X_u)|^2du\right)\\
&\leq C\left((N_t)^2+f(x)^2+\int_0^t F(X_u) du\right).
\end{align*}
Therefore,
\begin{align}\label{eq:supestimate}
\sup_{t \leq T} f(X_t)^2 \leq  C\left(\sup_{t \leq T} (N_t)^2+f(x)^2+\int_0^T F(X_u) du\right).
\end{align}
Due to~\eqref{eq:Lp}, $\sup_{t\leq T}|N_{t}| \in L^2$ and by~\eqref{eq:finite_moment}, we have $\mathbb{E}_x[\int_0^T F(X_u) du]< \infty$.
Hence we conclude that the right-hand side of~\eqref{eq:supestimate} is integrable. Summing over all $i$ yields integrability of $\sup_{t \leq T} F(X_t)$, which implies~\eqref{eq:sup} for all $f \in \mathcal{P}_m$.
\end{proof}

\begin{remark}\label{rem:remarkspol}
\begin{itemize}
\item[(i)] If $X$ is an $m$-polynomial process, then the process $Z=(X,X^{\mathbf{2}}, \ldots, X^{\mathbf{m}})$ is a $1$-polynomial process. If $m$ is even, the analysis of $m$-polynomial processes could be reduced to the study of $2$-polynomial processes at the cost of a more complicated state space, due to the construction $ Z'=(X,X^{\mathbf{2}}, \ldots, X^{\frac{\mathbf{m}}{2}})$.
\item[(ii)] Let us remark that the condition $m \geq 2$ in Theorem~\ref{th:characterizationpoly(iv)} is necessary for~\eqref{eq:extendedgenerator} being a true martingale. Indeed, the inverse $3$-dimensional Bessel process defined by $X= \frac{1}{\|B\|}$, where $B$ denotes a $3$-dimensional Brownian motion started at $B_0\neq 0$ satisfies
\[
dX_t=-X_t^2 dW_t, \quad X_0=x=\frac{1}{\|B_0\|},
\]
where $W$ is a one-dimensional standard Brownian motion. The extended generator is therefore given by 
\[
\mathcal{G}f(x)=\frac{1}{2}x^4\frac{d^2f(x)}{dx^2}.
\]
Hence $\mathcal{G}(\mathcal{P}_1)=0$, but
\begin{align}\label{eq:bessel}
X_t-x-\int_0^t\mathcal{G}X_s ds=X_t-x
\end{align}
is a strict local martingale, so $X$ is not a $1$-polynomial process.
\item[(iii)] In Definition~\ref{def-m-pol} we require $m$-polynomial processes to be also $k$-polynomial for all 
$k \in \{0, \ldots, m\}$, that is, we implicitly exclude processes whose extended generator maps polynomials of degree $k < m$ to polynomials of degree greater than $k \leq m$, while $\mathcal{G}(\PolmS) \subset \PolmS$ still holds true. Consider for instance 
\[
dX_t=\left(\frac{1}{2}-bX_t+\frac{1}{2}X_t^2\right)dt+\sqrt{X_t^2(1-X_t)}dW_t, \quad X_0=x \in [0,1],
\]
where $b \geq 1$ and $W$ is a one-dimensional standard Brownian motion. The state space is the interval $S=[0,1]$ and we have
\[
\mathcal{G}f(x)=\left(\frac{1}{2}-bx+\frac{1}{2}x^2\right)\frac{df(x)}{dx}+\frac{1}{2}x^2(1-x)\frac{d^2f(x)}{dx^2}.
\]
Thus $\mathcal{G}( \mathcal{P}_1)\subset  \mathcal{P}_2$, while $\mathcal{G}( \mathcal{P}_2)\subset  \mathcal{P}_2$.
Due to compactness of the state space it follows that $M^f$ is a true martingale for $f \in  \mathcal{P}_2$ and hence we can conclude that 
$P_t(\mathcal{P}_2)\subset  \mathcal{P}_2$ but $P_t( \mathcal{P}_1)\nsubseteq  \mathcal{P}_1$ as in the proof of Theorem~\ref{th:characterizationpoly} (iii) $\Rightarrow$ (i). On the other hand -- even though all moments of $ X $ exist and $ \mathcal{G}(\mathcal{P}_2) \subset \mathcal{P}_2 $ -- the subspace $ P_t (\mathcal{P}_3) $ is not only consisting of polynomials anymore. 
\end{itemize}
\end{remark}

\subsection{Polynomial processes and semimartingales}\label{sec:2.2}

The purpose of this section is to characterize polynomial processes as special semimartingales with 
characteristics of a particular form. Indeed, in Proposition~\ref{lem:semimartingale} below, we prove the equivalence of Condition~\ref{eq:cond_extended_gen} and the fact that the process
$(X_t1_{\{t < T_{\Delta}\}})$ is a special semimartingale whose characteristics are essentially
polynomials (of a particular degree) in $X$ and absolutely continuous with respect to the Lebesgue measure. 

Having established the particular form of the semimartingale characteristics, the characterization of the extended generator then simply follows from
It\^o's formula and is given by~\eqref{eq:real_generator}. This in turn allows us to state easy-to-verify conditions which guarantee that $X$ is an 
$m$-polynomial process for \emph{all} $m \geq 2$. Notice that this fills a gap which is left open in Theorem \ref{th:characterizationpoly(iv)}, where only the case of even $m$ is considered.

In order to formulate the following proposition concisely, we set
\[
 Y_t:=X_t1_{\{t < T_{\Delta}\}}=(f_1(X_t), \ldots, f_n(X_t))^{\top}
\]
and the define $C^2_m$ as the space of $2$-times continuously differentiable functions $g: S \to \mathbb{R}$, for which there exists some constant $\widetilde{C}$ such that
\begin{align*}
 |g(x)|+\sum_{i=1}^n |D_ig(x)|+\sum_{i,j=1}^n |D_{ij}g(x)|\leq \widetilde{C}(1+\|x\|^m).
\end{align*}

\begin{proposition}\label{lem:semimartingale}
Let $X$ be a time-homogeneous Markov process with state space $S_{\Delta}$ and let $m \geq 2$. Then the following assertions are equivalent.
\begin{itemize}
\item[(i)] Condition~\ref{eq:cond_extended_gen} holds, i.e.~$\PolmS \subset \mathcal{D}_{\mathcal{G}}$ and  $\mathcal{G}(\mathcal{P}_k) \subset \mathcal{P}_k$ for all $k \in \{0, 1, \ldots, m\}$.
 \item[(ii)]
$(Y_t)=(X_t1_{\{t< T_{\Delta}\}})$ is a semimartingale with respect to the stochastic basis $(\Omega, \mathcal{F}, (\mathcal{F}_t), \mathbb{P}_x)$. Moreover,
\begin{align}\label{eq:killing}
\mathbb{P}_x\left[t < T_{\Delta}\right]=e^{-\gamma t}
\end{align}
for some constant $\gamma \geq 0$ and the semimartingale characteristics $(B, C, \nu)$ associated with the ``truncation function'' $\chi(\xi) =\xi$ satisfy \footnote{All statements concerning the characteristics are meant up to an evanescent set.}
\begin{align}
 &B_{t,i}=\int_0^{t} b_i(X_{s}) ds,\label{eq:charb}\\
 & C_{t,ij}+\int_0^t\int_{\mathbb{R}^n}\xi_i\xi_j \nu(ds, d\xi)=\int_0^{t}  a_{ij}(X_{s}) ds, \label{eq:charcnu}
\end{align}
where $b_i \in \mathcal{P}_1$ and
$a_{ij} \in  \mathcal{P}_2$. Furthermore, the characteristics $C$ and $\nu$ can be written as
\begin{align}\label{eq:abscont}
C_{t,ij}=\int_{0}^{t}c_{ij}(X_s)ds, \quad \nu(\omega; dt, d\xi)=K(X_t(\omega),d\xi)dt,
\end{align}
where $c$ denotes some Borel measurable functions taking values in the set of positive semidefinite matrices 
and $K$ is a positive kernel from $(S, \mathcal{S})$ into $(\mathbb{R}^n, \mathcal{B}(\mathbb{R}^n))$ satisfying $K(x,\{0\})=0$ and $\int_{\{\|\xi\| \leq 1\}}\|\xi\|^2 K(x,d\xi)< \infty$.
Finally, we have for all $|\mathbf{k}| \in\{3, \ldots, m\}$
\begin{align}\label{eq:prop1}
\int_{\mathbb{R}^{n}} \xi^{\mathbf{k}} K(x,d\xi) =\sum_{|\mathbf{l}|=0}^{|\mathbf{k}|}\alpha_{\mathbf{l}} f^{\mathbf{l}}(x),
\end{align}
where $\alpha_{\mathbf{l}}$ denote some finite coefficients. 
\item[(iii)] $C_m^2$ lies in the domain of the extended generator of $X$ and for all $g \in C_m^2$, 
$\mathcal{G}$ is given by
\begin{equation}\label{eq:real_generator}
\begin{split}
 &\mathcal{G}g(x)=\sum_{i=1}^nD_i g(x)(b_{i}(x)+\gamma f_i(x))+\frac{1}{2}\sum_{i,j=1}^nD_{ij}g(x)c_{ij}(x)-\gamma g(x)\\
&\quad + \int_{\mathbb{R}^n} \left(g(x+\xi)-g(x)-\sum_{i=1}^nD_i g(x)\xi_i\right) (K(x,d\xi)-\gamma f_{\mathbf{0}}(x)\delta_{-x}(d\xi)),
\end{split}
\end{equation}
where $\gamma$ and $(b, c, K)$ satisfy the conditions of (ii) and $ f_{\mathbf{0}}(x)=1-1_{\Delta}(x)$.
\end{itemize}
All conditions (i), (ii) and (iii) imply that $Y$ is a \emph{special} semimartingale.
\end{proposition}

\begin{remark}\label{notation_pitfall}
\begin{itemize}
\item[(i)]
Concerning the direction (i) $\Rightarrow$ (ii), note that the existence of 
$\int\xi^{\mathbf{k}} K(x,d\xi)$ for all $|\mathbf{k}|\in\{2, \ldots, m\}$,
follows from the fact that $\mathcal{G}f(x)$ is a well defined polynomial for every $f \in \PolmS$.
In particular, it means that $\int\|\xi\|^kK(x,d\xi)< \infty$ for all $k\in\{2, \ldots, m\}$.
\item[(ii)]
Note that as a consequence of~\eqref{eq:charcnu} and~\eqref{eq:prop1} we have for all 
$k \in \{2, \ldots,  2\lfloor\frac{m}{2}\rfloor\}$ and $t \geq 0$ 
\begin{align}\label{eq:prop2}
\int_{\mathbb{R}^{n}} \|\xi\|^kK(X_t,d\xi) \leq \widetilde{C}\left(1+\|Y_t\|^{2\lfloor \frac{k+1}{2}\rfloor}\right),
\end{align}
with $\widetilde{C}$ some finite constant.
 \item[(iii)] 
Notice from the expression~\eqref{eq:real_generator} for $\mathcal{G}$ that 
the killing rate $\gamma$ is implicitly included in the compensator $K$ of the jump measure of $Y$ and also in $b$ due to the choice of the truncation function 
(compare~\cite[Section 3]{cfy}).
\end{itemize}
\end{remark}

\begin{proof}
Let us first prove $(i) \Rightarrow (ii)$.
Due to Condition~\ref{eq:cond_extended_gen} 
\begin{align}\label{eq:Mf}
M_t^f= f(X_t)-f(x)- \int_0^{t} \mathcal{G}f(X_s) ds,
\end{align}
is a $(\mathcal{F}_t,\mathbb{P}_x)$-local martingale for all $f \in \PolmS$. 
As the process
$
 \int_0^{t} \mathcal{G} f (X_s)ds
$
is predictable, $f(X)$ is a special $\mathbb{R}$-valued semimartingale for all $f \in \PolmS$. Note here that 
$
 f(X_t)\equiv f(X_t)1_{\{t < T_{\Delta}\}},
$
which is due to the convention $f(\Delta)=0$. In particular, $f(X_t)1_{\{t < T_{\Delta}\}}$ has c\`adl\`ag paths implying that $\lim_{t \rightarrow T_{\Delta}-}|f(X_t)|< \infty$ 
and $f(X)$ cannot explode. Choosing $f(x)=f_i(x)$ for $i=1, \ldots, n$, then implies that  
$(Y_t)=(X_t1_{\{t< T_{\Delta}\}})$ is an ($n$-dimensional) special semimartingale. 

Consider now~\eqref{eq:Mf} for $f=f_{\mathbf{0}}$. Since $\mathcal{G}f_{\mathbf{0}} \in \mathcal{P}_0$, $ M^{f_{\mathbf{0}}}$ is a true martingale and 
there exists some constant $\gamma$ such that
\[
\gamma 1_{\{t < T_{\Delta}\}}:=-\mathcal{G}f_{\mathbf{0}}(X_t).
\]
Taking expectations thus yields
\begin{align*}
\mathbb{E}_x\left[1_{\{t < T_{\Delta}\}}\right]&=\mathbb{E}_x\left[f_{\mathbf{0}}(X_t)\right]=1+\int_0^t \mathbb{E}_x\left[\mathcal{G}f_{\mathbf{0}}(X_s)\right]ds\\
&=1-\gamma \int_0^t \mathbb{E}_x\left[ 1_{\{s < T_{\Delta}\}}\right]ds,
\end{align*}
which in turn implies~\eqref{eq:killing}.

Let now $(B, C, \nu)$ denote the characteristics of $Y$ with respect to  the ``truncation function'' $\chi(\xi) =\xi$.
In order to determine their properties, we apply It\^o's formula
to $f_{\mathbf{k}}(X_t)$ for $k= |\mathbf{k}|\in \{1, \ldots, m\}$ for $ X_{t,i}=x_i+M^{f_i}_t+B_{t,i} $
\begin{equation}
\begin{split}\label{eq:Ito}
f_{\mathbf{k}}(X_t)&= f_{\mathbf{k}}(x)
+\int_{0}^{t}\sum_{i=1}^n D_i f_{\mathbf{k}}(X_{s-})dM^{f_i}_s
+\int_{0}^{t}\sum_{i=1}^nD_i f_{\mathbf{k}}(X_{s-})dB_{s,i}\\
&\quad +\frac{1}{2}\int_{0}^{t}\sum_{i,j=1}^nD_{ij}f_{\mathbf{k}}(X_{s-})dC_{s,ij} + \int_{0}^{t}\int_{\mathbb{R}^n} W(s, \xi)
\mu^{Y}(ds, d\xi),
\end{split}
\end{equation}
where  $\mu^{Y}$ denotes the random measure associated with the jumps of $Y$ and
\begin{align*}
 W(s,\xi):=\sum_{|\mathbf{l}|=2}^{k}\binom{\mathbf{k}}{\mathbf{l}}f_{\mathbf{k}-\mathbf{l}}(X_{s-})\xi^{\mathbf{l}}.
\end{align*}
Since $M^{f_i}$ is a local martingale and $X_{s-,i}$ is c\`agl\`ad,
$
(\int_{0}^{t}D_if_{\mathbf{k}}(X_{s-})dM^{f_i}_s)
$
is a local martingale, too, for all $ i \in \{1,\ldots,n\}$. Furthermore, the third and forth term on the right-hand side
are predictable processes of finite variation, thus in particular a process of locally integrable variation by~\cite[Lemma I.3.10]{jacod}. 
As $f_{\mathbf{k}}(X)$ is a special semimartingale, it follows from~\cite[Proposition I.4.22]{jacod} that
$
 \int_{0}^{t }\int_{\mathbb{R}^n} W(s, \xi) \mu^{Y}(ds, d\xi)
$
is also of locally integrable variation, since it is of finite variation and for a special semimartingale the finite variation part is locally integrable. 
Therefore, we have by~\cite[Proposition II.1.28]{jacod} that
\[
 \int_{0}^{t}\int_{\mathbb{R}^n} W(s, \xi) \mu^{Y}(ds, d\xi)-  \int_{0}^{t}\int_{\mathbb{R}^n} W(s, \xi) \nu(ds, d\xi)
\]
is a local martingale. Combining thus~\eqref{eq:Ito} with~\eqref{eq:Mf} and using the unique decomposition of a special semimartingale into a local martingale and a predictable finite variation process, we find
\begin{equation*}
\begin{split}
 M^{f_{\mathbf{k}}}&=\int_{0}^{t}\sum_{i=1}^n D_i f_{\mathbf{k}}(X_{s-})dM^{f_i}_s+
 \int_{0}^{t}\int_{\mathbb{R}^n} W(s, \xi) (\mu^{Y}(ds, d\xi)-\nu(ds,d\xi))\\
&=f_{\mathbf{k}}(X_t)-f_{\mathbf{k}}(x)
-\int_{0}^{t}\sum_{i=1}^nD_i f_{\mathbf{k}}(X_{s-})dB_{s,i}\\
&\quad-\int_{0}^{t}\frac{1}{2}\sum_{i,j=1}^nD_{ij}f_{\mathbf{k}}(X_{s-})dC_{s,ij} - \int_{0}^{t}\int_{\mathbb{R}^n} W(s, \xi)
\nu(ds, d\xi).
\end{split}
\end{equation*}
Therefore 
\begin{equation}
\begin{split}\label{eq:gendiffop}
\int_0^{t} \mathcal{G}f_{\mathbf{k}}(X_s)ds&=\int_{0}^{t}\sum_{i=1}^nD_i f_{\mathbf{k}}(X_{s})dB_{s,i}+\int_{0}^{t}\frac{1}{2}\sum_{i,j=1}^nD_{ij}f_{\mathbf{k}}(X_{s-})dC_{s,ij}\\
&\quad + \int_{0}^{t}\int_{\mathbb{R}^n} W(s, \xi)\nu(ds, d\xi).
\end{split}
\end{equation}
Consider now~\eqref{eq:gendiffop} for $|\mathbf{k}|=1$, i.e.~the polynomials $f_{i}(x)$, where $i\in \{1, \ldots, n\}$.
In this case~\eqref{eq:gendiffop} reads as
\[
\int_0^{t} \mathcal{G}f_i(X_{s})ds=B_{t,i}.
\]
Setting $b_i(x):=\mathcal{G}f_i(x)$ 
therefore implies that $b_i \in  \mathcal{P}_1$.
Moreover, applying ~\eqref{eq:gendiffop} to the quadratic polynomials $f(x)=f_{ij}(x)$ for $i,j\in \{1, \ldots, n\}$
yields
\begin{align}\label{eq:Cnu}
 C_{t,ij}+\int_{0}^{t}\int_{\mathbb{R}^n} \xi_i\xi_j\nu(ds, d\xi)=\int_0^{t} a_{ij}(X_s) ds
\end{align}
for some $a_{ij}\in  \mathcal{P}_2$, since $\mathcal{G}f_{ij}(x)$ and $D_if_{ij}(x)b_i(x)=f_j(x)b_i(x)$ lie in $ \mathcal{P}_2$. Hence we have proved~\eqref{eq:charb} and~\eqref{eq:charcnu}.

In order to show~\eqref{eq:abscont}, we define $A'_t(\omega)=\int_{0}^{t}\int_{\mathbb{R}^n} \|\xi\|^2\nu(\omega; ds, d\xi)$. By the same arguments as in the proof of~\cite[Proposition II.2.9.b]{jacod}, there exists a random measure $K'(\omega, t; d\xi)$ on $(\mathbb{R}^n, \mathcal{B}(\mathbb{R}^n))$ such that $\nu(\omega; dt, d\xi)=K'(\omega, t; d\xi)dA_t(\omega)$. Moreover, since
\[
\sum_{i=1}^n C_{t,ii}(\omega)+A'_t(\omega)=\sum_{i=1}^n\int_0^{t} a_{ii}(X_s(\omega)) ds=:\int_0^{t}a_s(\omega) ds
\]
and as $C_{t,ii}, \, i \in \{1,\ldots,n\}$, and $A'_t$ are non-negative increasing processes (of finite variation), 
$C_{ii}$ and $A'$ are absolutely continuous with respect to Lebesgue measure. 
Hence,~\cite[Proposition I.3.13]{jacod} implies the existence of predictable processes $\widetilde{c}_{ii}$ and $H$ such that 
$C_{t,ii}=\int_0^t\widetilde{c}_{s,ii}ds$ and $A'_t=\int_0^tH_s ds$. Then $\widetilde{K}_{\omega,t}(d\xi)=H_t(\omega)K'(\omega,t; d\xi)$ 
is again a predictable random measure satisfying $\nu(\omega; dt,d\xi)=\widetilde{K}_{\omega,t}(d\xi)dt$ almost surely. 
Having constructed this kernel,~\eqref{eq:Cnu} now becomes
\[
 C_{t,ij}=\int_0^{t} \left(a_{ij}(X_s)-\int_{\mathbb{R}^n} \xi_i\xi_j \widetilde{K}_{\omega,s}(d\xi)\right)ds,
\]
implying that $C_{t,ij}$ for $i\neq j$ is also absolutely continuous with respect to Lebesgue measure and can therefore be written as $C_{t,ij}=\int_0^t\widetilde{c}_{s,ij}ds$.
Finally, by~\cite[Theorem 6.27]{cinlar} we can choose homogeneous versions for the processes $\widetilde{c}$ and $\widetilde{K}$ such that
$
 \widetilde{c}_t(\omega)=c(X_t(\omega)) \textrm{ and }  \widetilde{K}_{\omega,t}(d\xi)=K(X_t(\omega),d\xi).
$

It remains to establish property~\eqref{eq:prop1}. To this end notice that~\eqref{eq:gendiffop} can be written as
\begin{align*}
\mathcal{G}f_{\mathbf{k}}(x)&=\sum_{i=1}^nD_i f_{\mathbf{k}}(x)b_i(x)+\frac{1}{2}\sum_{i,j=1}^nD_{ij}f_{\mathbf{k}}(x)\left(c_{ij}(x)+\int_{\mathbb{R}^n}\xi_i \xi_j K(x,d\xi)\right)\\
&\quad+\int_{\mathbb{R}^n} \left(\sum_{|\mathbf{l}|=3}^{k}\binom{\mathbf{k}}{\mathbf{l}}f_{\mathbf{k}-\mathbf{l}}(x)\xi^{\mathbf{l}}\right)K(x,d\xi).
\end{align*}
Since $\mathcal{G}f_{\mathbf{k}}(x)$, $ D_if_{\mathbf{k}}(x)b(x)$ and $D_{ij}f_{\mathbf{k}}(x)(c_{ij}(x)+\int_{\mathbb{R}^n}\xi_i \xi_jK(x,d\xi))$ 
lie in $\mathcal{P}_k$ for all $k=|\mathbf{k}|\leq m$,~\eqref{eq:prop1} simply follows by induction. 

Let us now prove the implication (ii) $\Rightarrow$ (iii). For notational simplicity we set
\[
 V(Y_s,\xi):=g(Y_s+\xi)-g(Y_s)-\sum_{i=1}^nD_i g(Y_s)\xi_i
\]
for $g \in C_m^2$. Then, since $\int \|\xi\|^k  K(x,d\xi)< \infty $ for all $k \in \{2, \ldots,m\}$, which is a consequence of~\eqref{eq:charcnu}  and~\eqref{eq:prop1},
we have 
\[
  \int_{\mathbb{R}^n} \left|V(Y_s,\xi)\right| K(X_s,d\xi)
\leq h(Y_s)+H(Y_s)\int_{\mathbb{R}^n}  (\|\xi\|^2 \wedge \|\xi\|^m) K(X_s,d\xi)< \infty, 
\]
where $h$ and $H$ denote some positive (finite-valued) functions. 
Hence, the process $\int_0^{\cdot} \int V(Y_s,d\xi) K(X_s,d\xi)ds$ is of locally integrable variation
and It\^o's formula thus implies that
\begin{equation}\label{eq:mathcalM}
\begin{split}
\mathcal{M}^g_t&:=g(Y_t)-g(x)
-\int_{0}^{t}\sum_{i=1}^nD_i g(Y_{s})b_{i}(X_s)ds-\int_{0}^{t}\frac{1}{2}\sum_{i,j=1}^nD_{ij}g(Y_{s})c_{ij}(X_s)\\
&\qquad - \int_{0}^{t}\int_{\mathbb{R}^n}   V(Y_s,\xi) K(X_s,d\xi)ds
\end{split}
\end{equation}
is a local martingale (compare also~\cite[Theorem II.2.42]{jacod}).
Moreover, 
\[
-1_{\{ T_{\Delta}\leq t\}}+\gamma(t \wedge T_{\Delta} ) =1_{\{t < T_{\Delta}\}}-1 +\int_0^t \gamma 1_{\{s < T_{\Delta}\}} ds
\]
is a martingale, since
\begin{align*}
\mathbb{E}_{x}\left[1_{\{t < T_{\Delta}\}}-1 +\int_0^t \gamma 1_{\{s < T_{\Delta}\}} ds \right]
&=e^{-\gamma t}-1+\int_0^{t} \gamma e^{-\gamma s} ds=0.
\end{align*}
Denote now by $\mathcal{G}^{\#}$ the right-hand side of~\eqref{eq:real_generator}. We need to prove that
\[
M_t^{\#g}:= g(X_t)-g(x)-\int_0^t \mathcal{G}^{\#}g(X_s)ds 
\]
is a local martingale. By the definition of $\mathcal{G}^{\#}$ we have
\[
 M_t^{\#g}=\mathcal{M}_t^g-g(0)\left(1_{\{T_{\Delta}\leq t\}}- \gamma(t \wedge T_{\Delta} )\right),
\]
where $\mathcal{M}^g$ is given by~\eqref{eq:mathcalM}. Since both terms on the right-hand side are local martingales, the same holds true for  $M^{\#g}$.

Finally (i) follows from (iii), since $\mathcal{P}_m \subset C_m^2$ and since $\mathcal{G}$ applied to $f_{\mathbf{k}}$ 
maps $\mathcal{P}_k$ into $\mathcal{P}_k$ for $ k\in \{0, \ldots, m\}$, which is due to the assumptions on the characteristics.
\end{proof}

Since every $m$-polynomial process satisfies Condition~\ref{eq:cond_extended_gen}, the following corollary is an obvious consequence of Proposition~\ref{lem:semimartingale}.

\begin{corollary}
 Let $X$ be an $m$-polynomial process with $m \geq 2$. 
Then $(Y_t)=(X_t1_{\{t < T_{\Delta}\}})$ is a special semimartingale satisfying the conditions~\eqref{eq:killing} -~\eqref{eq:prop1} and its extended generator is of 
form~\eqref{eq:real_generator}.
\end{corollary}

When $m$ is an even number, the converse direction also follows easily from Theorem~\ref{th:characterizationpoly(iv)} and Proposition~\ref{lem:semimartingale}.
The general case is treated in the subsequent theorem, where we  provide sufficient conditions in terms of the compensator of the jump measure such that the converse statement 
holds true. Its proof relies on a maximal inequality which can be established for semimartingales whose characteristics satisfy 
the conditions~\eqref{eq:charb} -~\eqref{eq:prop1}. This is subject of Lemma~\ref{lem:supX} below.

\begin{theorem}\label{cor:suffcond}
Let $X$ be a time-homogeneous Markov process with state space $S_{\Delta}$ and let $m \geq 2$. Suppose that $(Y_t)=(X_t1_{\{t <T_{\Delta}\}})$ is a semimartingale,
which satisfies the conditions~\eqref{eq:killing} -~\eqref{eq:prop1} (with respect to the ``truncation function'' $\chi(\xi)=\xi$) 
or equivalently that $C^2_m \subset \mathcal{D}_{\mathcal{G}}$ and that its extended generator $\mathcal{G}$ is given by~\eqref{eq:real_generator}.
If 
\begin{align}\label{eq:mmomentjump}
 \mathbb{E}_x\left[\int_{\mathbb{R}^n}\|\xi\|^m K(X_t,d\xi)\right]< \infty, \quad \textrm{for almost all }  t \geq 0,
\end{align}
or if 
\begin{align}\label{eq:Xdomjump}
 \int_{\mathbb{R}^{n}} \|\xi\|^mK(X_t,d\xi) \leq \widetilde{C}\left(1+\|Y_t\|^{m}\right),  \quad t \geq 0
\end{align}
for some constant $\widetilde{C}$, then $X$ is an $m$-polynomial process.
\end{theorem}

\begin{proof}
Similarly as in the proof of Theorem~\ref{th:characterizationpoly(iv)}, it suffices to show that
for each $f$, $x \in S_{\Delta}$ and every fixed $t \geq 0$, $P_t|f|(x)=\mathbb{E}_x[|f(X_t)|]< \infty$ and
\begin{align*}
\mathbb{E}_x\left[\sup_{s\leq t}\left|M_s^{f}\right|\right]< \infty.
\end{align*}
Due to the assumptions~\eqref{eq:mmomentjump} or~\eqref{eq:Xdomjump}, this follows from the moment estimate~\eqref{eq:momentestimate} proved in Lemma~\ref{lem:supX} below.
\end{proof}

\begin{remark}
Let us remark that under~\eqref{eq:prop1}, condition~\eqref{eq:Xdomjump} is always satisfied when $m$ is an even number (see also Remark~\ref{notation_pitfall} (ii)).
In this case 
$
\mathbb{E}_x\left[\sup_{s\leq t}\left|M_s^f\right|\right]
$
is always finite, as already shown in the proof of Theorem~\ref{th:characterizationpoly(iv)}.
If $m > 2$ is an odd number, Condition~\ref{eq:cond_extended_gen} together with~\eqref{eq:mmomentjump} or~\eqref{eq:Xdomjump}
is sufficient for $X$ being an $m$-polynomial process.
\end{remark}

Using the structure of the semimartingale characteristics derived in Proposition~\ref{lem:semimartingale}, we now state the announced maximal inequality 
and some associated moment estimates. This result is probably known but the proof is included for convenience. 
A similar statement for the case of L\'evy driven SDEs can be found in~\cite{jacodkurtz}.

\begin{lemma}\label{lem:supX}
Fix $t>0$ and let $m \geq 2$. Let $Y$ be a semimartingale with respect to $(\Omega, \mathcal{F}, (\mathcal{F}_t), \mathbb{P}_x)$, 
whose characteristics $(B,C, \nu)$ associated with the ``truncation function'' $\chi(\xi)=\xi$ satisfy the conditions~\eqref{eq:charb},~\eqref{eq:charcnu} and~\eqref{eq:abscont} given in
Proposition~\ref{lem:semimartingale}. 
Then there exists a constant $\widetilde{C}$ such that
\begin{equation}
\begin{split} \label{eq:maxinequality}
\mathbb{E}_x\left[\sup_{s\leq t}\left\|Y_s\right\|^m\right]&\leq 
\widetilde{C}\Bigg(\|x\|^m+1+\int_0^t \mathbb{E}_x\left[\int_{\mathbb{R}^n}\|\xi\|^m K(X_s,d\xi)\right]ds\\
&\quad \quad+\int_0^t\mathbb{E}_x\left[\|Y_{s}\|^m\right]ds\Bigg).
\end{split}
\end{equation}
In particular, if one of the conditions~\eqref{eq:mmomentjump} or~\eqref{eq:Xdomjump} is satisfied, then
there exist finite constants $K$ and $\widetilde{C}$ such that 
\begin{align}\label{eq:momentestimate}
\mathbb{E}_x\left[\sup_{s\leq t}\left\|Y_s\right\|^m\right] \leq Ke^{\widetilde{C}t}.
\end{align}
\end{lemma}

\begin{proof}
For notational simplicity we only consider the case where $Y$ is one-dimensional and thus omit all indices.
Due to Proposition~\ref{lem:semimartingale} and the assumptions on the characteristics, $Y$ is a special semimartingale and its 
canonical decomposition is given by $Y=x+M+\int_0^{\cdot} b(X_u)du$, where we write $M$ for $M^{f_1}$.
Denote by $Z$ the quadratic variation of the purely discontinuous martingale part of $Y$, that is, 
\[
 Z_{t}=\sum_{s \leq t} (\Delta Y_{s})^2=\int_0^t \int_{\mathbb{R}} \xi^2 \mu^Y(ds,d\xi),
\]
where $\mu^Y$ is the random measure associated with the jumps of $Y$. Define furthermore stopping times
\begin{align*}
 T_j^{Y}&=\inf\{t \geq 0 \,|\, |Y_{t}| \geq j \textrm{ or } |Y_{t-}| \geq j \}\\
 T_j^Z&=\inf\{t \geq 0 \,|\, |Z_{t}| \geq j \textrm{ or } |Z_{t-}| \geq j \}
\end{align*}
and set $T_j =T_j^{Y} \wedge T_j^X$. We can estimate
\[
\sup_{s\leq t}\left|Y^{T_j}_{s}\right|^m\leq \widetilde{C}\left(|x|^m+\sup_{s\leq t}\left|M^{T_j}_{s}\right|^m+
\sup_{s\leq t}\left|\int_0^{s \wedge {T_j}} b(X_u)du\right|^m\right),
\]
where $\widetilde{C}$ denotes some constant which may vary from line to line.

Since $b \in  \mathcal{P}_1$ we have
\begin{equation}\label{eq:estimateb}
\begin{split}
\sup_{s\leq t}\left|\int_0^{s  \wedge {T_j}}b(X_{u}) du\right|^m
&\leq \widetilde{C}\left(1+\int_0^{t  }|Y_u^{T_j}|^m du \right).
\end{split}
\end{equation}
Concerning $\sup_{s \leq t}|M^{T_j}_{s}|$, an application of Burkholder-Davis-Gundy's inequality yields
\begin{align}\label{eq:estimateC}
\mathbb{E}_x\left[\sup_{s\leq t}\left|M^{T_j}_{s}\right|^m\right]\leq \widetilde{C}\mathbb{E}_x\left[\left[M,M\right]_{t\wedge{T_j}}^{\frac{m}{2}}\right]
\leq \widetilde{C}\mathbb{E}_x\left[C_{t\wedge{T_j}}^{\frac{m}{2}}+Z_{t\wedge{T_j}}^{\frac{m}{2}}\right]
\end{align}
As $C$ satisfies~\eqref{eq:charcnu}, we can estimate it by $C_{t} \leq \int_0^{t }a(X_s)ds$,
where $a \in  \mathcal{P}_2$ is nonnegative, and we get
\[
\mathbb{E}_x\left[C_{t\wedge{T_j}}^{\frac{m}{2}} \right]\leq \widetilde{C}\left( 1+\int_0^t \mathbb{E}_x\left[|Y_s^{T_j}|^m\right]ds\right).
\]
Therefore it remains to handle $
Z_{t\wedge{T_j}}^{\frac{m}{2}}$. Following the approach of~\cite{jacodkurtz}, we can write
\begin{align*}
Z_{t\wedge T_j}^{\frac{m}{2}}&=\sum_{s\leq t \wedge T_j}(Z_{s_-}+\Delta Z_s)^{\frac{m}{2}}-(Z_{s_-})^{\frac{m}{2}}\\
&=\int_0^{t \wedge T_j}\int_{\re^n}\left((Z_{s_-}+\xi^2)^{\frac{m}{2}}-(Z_{s_-})^{\frac{m}{2}}\right)\mu^{Y}(ds,d\xi),
\end{align*}
which is due to the fact that $Z$ is purely discontinuous, non-decreasing and 
$\Delta Z_s=|\Delta Y_s|^2$. Furthermore, since $\nu$ is the predictable compensator of $\mu^{Y}$, we have
\begin{align}\label{eq:Yexpectation}
\mathbb{E}_x\left[Z_{t\wedge T_j}^{\frac{m}{2}}\right]=
\mathbb{E}_x\left[\int_0^{t \wedge T_j}\int_{\re^n}\left((Z_{s_-}+\xi^2)^{\frac{m}{2}}-(Z_{s_-})^{\frac{m}{2}}\right)\nu(ds,d\xi)\right].
\end{align}
In the sequel we shall use the following inequalities (see~\cite{jacodkurtz})
\begin{align}
(z+x)^p-z^p&\leq2^{p-1}(z^{p-1}x+x^p)\label{eq:ineq1},\\
z^{p-1}x&\leq \varepsilon z^p+\frac{x^p}{\varepsilon^{p-1}},\label{eq:ineq2}
\end{align}
for $x,z \geq 0$, $\varepsilon >0$ and $p \geq 1$.
Applying~\eqref{eq:ineq1}, equation~\eqref{eq:Yexpectation} becomes
\[
\mathbb{E}_x\left[Z_{t\wedge T_j}^{\frac{m}{2}}\right] \le 
\mathbb{E}_x\left[\int_0^{t \wedge T_j}\int_{\re^n}2^{\frac{m}{2}-1}\left(Z_{s-}^{\frac{m}{2}-1}\xi^2+|\xi|^m\right)\nu(ds,d\xi)\right].
\]
For the first part, we then have due to the assumption on $\nu$
\begin{multline*}
\mathbb{E}_x\left[\int_0^{t \wedge T_j}\int_{\re^n} 2^{\frac{m}{2}-1}Z_{s-}^{\frac{m}{2}-1}\xi^2 \nu(ds,d\xi)\right]\\
= \mathbb{E}_x\left[\int_0^{t \wedge T_j}2^{\frac{m}{2}-1}Z_{s}^{\frac{m}{2}-1}\left(\int_{\re^n} \xi^2 K(X_s,d\xi)\right) ds\right]\\
\leq  \mathbb{E}_x\left[\int_0^{t \wedge T_j} 2^{\frac{m}{2}-1}Z_{s}^{\frac{m}{2}-1}a(X_s) ds\right]\\
\leq \mathbb{E}_x\left[\int_0^{t \wedge T_j } \widetilde{C} \left(\varepsilon Z_{s}^{\frac{m}{2}}+
\frac{1+|Y_s|^m}{\varepsilon^{\frac{m}{2}-1}}\right)ds\right],
\end{multline*}
where we use the fact that $Z$ is nonnegative, $a \in \mathcal{P}_2$ and~\eqref{eq:ineq2}.
Estimating
 $
 \int_0^{t \wedge T_j }Z_{s}^{\frac{m}{2}}ds\leq  j^{\frac{m}{2}} \wedge Z_{t\wedge T_j}^{\frac{m}{2}},
 $
which follows from the fact that $Z_{s}$ is non-decreasing and $Z_{s-} \leq j$ for $s \leq T_j$, we finally obtain
\begin{align*}
 \mathbb{E}_x\left[Z_{t\wedge T_j}^{\frac{m}{2}}\right]
&\leq \widetilde{C} \varepsilon \mathbb{E}_x\left[j^{\frac{m}{2}} \wedge Z_{t \wedge T_n}^{\frac{m}{2}}\right]
+\mathbb{E}_x\left[\int_0^{t \wedge T_j } \frac{\widetilde{C}} {\varepsilon^{\frac{m}{2}-1}}(1+|Y_s|^m)ds\right]\\
&\quad+\mathbb{E}_x\left[\int_0^{t \wedge T_j}\int_{\re^n}2^{\frac{m}{2}-1}|\xi|^m K(X_s,d\xi)ds\right].
\end{align*}
 Choosing $\varepsilon=\frac{1}{2\widetilde{C}}$ leads to
 \begin{align*}
 \frac{1}{2}\mathbb{E}_x\left[Z_{t \wedge T_j}^{\frac{m}{2}}\right]&\leq\mathbb{E}_x\left[Z_{t\wedge T_j}^{\frac{m}{2}}\right]-\frac{1}{2}\mathbb{E}_x\left[j^{\frac{m}{2}} \wedge Z_{t\wedge T_j}^{\frac{m}{2}}\right] \\
 &\leq \widetilde{C}\mathbb{E}_x\left[\int_0^{t}\left(1+|Y^{T_j}_{s}|^m+\int_{\re^n}|\xi|^m K(X_{s \wedge T_j},d\xi)\right)ds\right].
 \end{align*}
Combining this with the estimates~\eqref{eq:estimateb} and~\eqref{eq:estimateC}, we find
\begin{equation}\label{eq:stoppedmaxineq}
 \begin{split}
\mathbb{E}_x\left[\sup_{s\leq t}\left|Y^{T_j}_{s}\right|^m\right] 
&\leq  \widetilde{C}\Bigg(|x|^m+1+\int_0^t\mathbb{E}_x\left[\int_{\mathbb{R}}|\xi|^m K(X_{s \wedge T_j},d\xi)\right]ds\\
&\quad \quad \int_0^t \mathbb{E}_x\left[|Y^{T_j}_{s}|^m\right]ds\Bigg).
\end{split}
\end{equation}
By monotone convergence we obtain~\eqref{eq:maxinequality}.

Concerning~\eqref{eq:momentestimate}, note that under the conditions~\eqref{eq:mmomentjump} or~\eqref{eq:Xdomjump}, we can deduce from~\eqref{eq:stoppedmaxineq} that
\begin{align*}
\mathbb{E}_x\left[\sup_{s\leq t}\left|Y^{T_j}_{s}\right|^m\right] \leq K+ \widetilde{C}\int_0^t \mathbb{E}_x\left[\sup_{u \leq s}|Y^{T_j}_{u}|^m\right]ds
\end{align*}
for some finite constants $K$ and $\widetilde{C}$.
Since the right hand side is finite due to the conditions~\eqref{eq:mmomentjump} or~\eqref{eq:Xdomjump} and the estimate
\begin{align*}
 \sup_{u \leq s}|Y^{T_j}_{u}|^m&\leq |x|^m +j^m + \mathbb{E}_x\left[|\Delta Y_{s \wedge T_j}|^m\right]\\
&\leq |x|^m +j^m + \mathbb{E}_x\left[\int_0^{s \wedge T_j} \int_{\mathbb{R}}|\xi|^m K(X_{s},d\xi) ds\right]< \infty,
\end{align*}
Gronwall's lemma yields
\[
 \mathbb{E}_x\left[\sup_{s\leq t}\left|Y^{T_j}_{s}\right|^m\right]\leq Ke^{\widetilde{C}t}
\]
for all $j \in \mathbb{N}$ and the result follows from monotone convergence.
\end{proof}

\begin{remark}
\begin{itemize}
\item[(i)]
It is important to note that the characteristics of $(Y)=(X_t1_{\{t <T_{\Delta}\}})$ in the above statements are always specified with respect to the 
``truncation function'' $\chi(\xi)=\xi$. While $C$ and $\nu$ do not depend on this choice, the characteristic $ B $ does depend on $ \chi $.  
So, if one chooses another truncation function $\chi'$, then the difference between $ B $ and $ B'$ is given by
$\int_0^t\int_{\re^n\setminus\{0\}}\left(\chi'(\xi)-\chi(\xi)\right)\nu(ds,d\xi)$.
Thus, the requirement that $C$ and $\nu$ are as in Theorem~\ref{cor:suffcond} and
 \begin{align}\label{eq:truncation}
  \left(b(X_t)+\int_{\re^n\setminus\{0\}}\left(\chi(\xi)-\chi'(\xi)\right)K(X_t,d\xi)\right)=\sum_{|\mathbf{k}|=0}^1\alpha_{\mathbf{k}}f_{\mathbf{k}}(X_t)
\end{align}
is an equivalent condition guaranteeing that $X$ is $m$-polynomial.
\item[(ii)]\label{rem:example}
We now give two examples of kernels $K(x,d\xi)$ which satisfy the conditions of Proposition~\ref{lem:semimartingale} as long as 
$c_{ij} \in  \mathcal{P}_2$ and $b$ satisfies~\eqref{eq:truncation}.

\begin{itemize}
\item[(a)] The first one essentially requires $K(x, d\xi)$ to be a quadratic polynomial in $x$, that is,
\begin{align*}
K(x, d\xi)=\left(\frac{\mu_{00}(d\xi)}{\|\xi\|^2 \wedge 1}+\sum_{i=1}^{n}x_{i}\frac{\mu_{i0}(d\xi)}{\|\xi\|^2 \wedge 1}+\sum_{i\leq j}x_{i}x_{j}\frac{\mu_{ij}(d\xi)}{\|\xi\|^2 \wedge 1}\right),
\end{align*}
where all $\mu_{ij}$ are finite signed measures on $\re^n$ such that $K(x,\cdot)$ is a well defined L\'evy measure for every $x \in S$.
In view of Remark~\ref{notation_pitfall}, it is necessary to require
\begin{align*}
\int_{\|\xi\|>1} \|\xi\|^m (\mu^+_{ij}(d\xi)+\mu_{ij}^-(d\xi))=\int_{\|\xi\|>1} \|\xi\|^m  |\mu_{ij}(d\xi)|< \infty,
\end{align*}
where $\mu_{ij}^+, \mu_{ij}^{-}$ denotes the Jordan decomposition of $\mu_{ij}$.
\item[(b)] Alternatively, $K$ can be specified as the pushforward of a L\'evy measure under an affine function.
Let $d \geq 1$ and let
\[
p: S \times \re^d \rightarrow \re^n, \quad (x,y) \mapsto p(x,y)=p^x(y)=H(y)x+h(y),
\]
be an affine function in $x$. Here, $H: \re^d \rightarrow \re^{n\times n}$ and $h: \re^d \rightarrow \re^n$ are assumed to be measurable.
We define $K$ then by
\begin{align*}
 K(x,d\xi):={(p^{x})}_{\ast}\mu(d\xi),
\end{align*}
where for each $x \in S$, ${(p^x)}_{\ast}\mu$ denotes the pushforward of the measure $\mu$ under the map $p^x$. Moreover, $\mu$ is a L\'evy measure on $\re^d$ integrating
\begin{align*}
\int_{\re^d \setminus\{0\}}\left(\|H(y)\|^k+\|h(y)\|^k\right) \mu(dy)  \quad \textrm{ for all } k \in \{1, \ldots, m\}.
\end{align*}
\end{itemize}
\end{itemize}
\end{remark}

\section{Examples}\label{sec:examples}

In order to apply Theorem~\ref{cor:suffcond} to the following examples, we assume $m \geq 2$.

\begin{example}[Affine processes]\label{ex:affineprocesses}
Every affine process $X$ on $S=\mathbb{R}_{+}^p \times \mathbb{R}^{n-p} $ 
is $m$-polynomial if the killing rate is constant and if the L\'evy measures $\mu_i$ for $i \in \{0,1,\ldots,p\}$, satisfy
\begin{align}\label{eq:affinejumps}
\int_{{\|\xi\|>1}} \|\xi\|^m \mu_i(d\xi)< \infty.
\end{align}
For details on affine processes see, e.g.,~\cite{dfs} and~\cite{kst}.\footnote{We write here $\mu_0$ for the constant part of the jump measure in contrast to~\cite{dfs}, where it is denoted by $m$.}

\begin{proof}
Since the differential characteristics of $(Y)=(X_t1_{\{t < T_{\Delta}\}})$ are affine functions in $X$ and since~\eqref{eq:affinejumps} assures that
condition~\eqref{eq:Xdomjump} is satisfied, that is,
\begin{align*}
 \int_{\mathbb{R}^n} \|\xi\|^m K(X_t,d\xi)&= \int_{\mathbb{R}^n}\|\xi\|^m \left(\mu_0(d\xi)+\sum_{i=1}^p X_{t,i}\mu_i(d\xi)\right)\\
&\leq \widetilde{C} (1+ \|Y_t\|) \leq \widetilde{C} (1+ \|Y_t\|^m), 
\end{align*}
for some constant $\widetilde{C}$ and $m \geq 2$, the assertion follows from Theorem~\ref{cor:suffcond}.
 \end{proof}

Let us remark that affine processes are defined via their characteristic function, which is of the form
\[
 \mathbb{E}_x\left[e^{\langle \im  u, X_t\rangle}\right]=e^{\phi(t,\im u)+\langle\psi(t,\im u),x\rangle}
\]
for some function $\phi$ and $\psi$.
This definition then implies affine semimartingale characteristics, from which the polynomial property can easily be seen due to Theorem~\ref{cor:suffcond}.
Note also that the explicit knowledge of $\phi$ and $\psi$ is not necessary to compute the moments of an affine process. 
Simply the knowledge of its characteristics, which determine the linear map $A$ as given in Theorem~\ref{th:characterizationpoly} (ii), is enough (see Section~\ref{sec-applications}).
Moreover, for affine processes, it has not been proved so far that the existence of the moments of the L\'evy measures automatically implies the existence of the moments 
of the process itself. We can conclude this simply from the more general statement of Lemma~\ref{lem:supX}.
\end{example}

\begin{example}[Exponential L\'evy models]
Exponential L\'evy models are of the form
$X=xe^{L}$, where $L$ is a L\'evy process on $\re$ with triplet $(b,c,\mu)$. Under the integrability assumption
$\int_{|y|>1}e^{my}\mu(dy) < \infty$,
which guarantees the existence of $\mathbb{E}_x\left[|X_t|^m\right]$,
exponential L\'evy models are $m$-polynomial, since we have
$
\mathbb{E}_x\left[x^me^{mL_t}\right]=x^me^{t\psi(m)},
$
where $\psi$ denotes the cumulant generating function of the L\'evy process.
\end{example}

\begin{example}[L\'evy driven SDEs]
Let $L$ denote a L\'evy process on $\re^d$ with triplet $(b,c,\mu)$. Suppose furthermore that $ V_1,\ldots,V_d $ are affine functions, i.e.~we have
$
V_i: S\rightarrow \re^n,\, x\mapsto H_ix+h_i,
$
where $H_i \in \re^{n\times n}$ and $h_i \in \re^n$. A process $X$ which solves the stochastic differential equation
\begin{align*}
dX_t = \sum_{i=1}^d V_i (X_{t-}) dL_{t,i}, \quad X_0 = x \in S,
\end{align*}
and leaves $S$ invariant, is $m$-polynomial, if the following moment condition on the L\'evy measure 
\begin{align}\label{eq:LevySDE}
\int_{{\|y\|>1}} \|y\|^{m} \mu(dy)< \infty
\end{align}
is satisfied.

\begin{proof}
For $C_m^2$-functions $g$ and general Lipschitz continuous functions $V_1,\ldots,V_d$ the extended generator of $X$ with respect to some truncation function $\chi$ is given by
\begin{align*}
&\left\langle\sum_{i=1}^d V_i(x) b_i,\nabla g(x)\right\rangle+ \frac{1}{2}\sum_{i,j=1}^n ((V_1(x)\ldots V_d(x))c(V_1(x)\ldots V_d(x))')_{ij}D_{ij}g(x)\\
&+\int\left(g\left(x+\sum_{i=1}^d V_i(x)y_i\right)-g(x)-\left\langle\sum_{i=1}^{d}V_i(x)\chi_i(y),\nabla g(x)\right\rangle\right) \mu(dy).
\end{align*}
Concerning the the compensator of the jump measure $K(x,d\xi)$, this example corresponds to the situation of Remark~\ref{rem:example} (ii) (b) 
with $p(x,y)=H(y)x+h(y)=\sum_{i=1}^d H_iy_ix+h_iy_i$. 
Condition~\eqref{eq:Xdomjump}, that is,
\[
 \int_{\mathbb{R}^n} \|\xi\|^m K(X_t,d\xi)=\int_{\mathbb{R}^n} \left\|\sum_{i=1}^d H_iy_iX_t+h_iy_i\right\|^m \mu(dy)\leq \widetilde{C} (1+ \|X_t\|^m)
\]
for some constant $\widetilde{C}$, is satisfied due to~\eqref{eq:LevySDE}. Hence Theorem~\ref{cor:suffcond} yields the assertion.
\end{proof}
\end{example}

\begin{example}[Quadratic term structure models~\cite{chen}]
Consider the following quadratic term structure model $r$, specified as non-negative quadratic function of a one-dimensional Ornstein-Uhlenbeck process $Y$
\[
r_t=R_0+R_1Y_t+R_2Y_t^2,
\]
for appropriate $R_i \in \re$. Here, $Y$ is given by
\[
dY_t=(b+\beta Y_t)dt+\sigma dW_t,
\]
where $W$ is a standard Brownian motion. The joint process $X=(Y, r)$ then satisfies the dynamics
\begin{align*}
\left( \begin{array}{c}
dY_t\\
dr_t\end{array}\right)&=
\left(\left( \begin{array}{c}
b\\
R_1b+R_2\sigma^2-2R_0\beta
\end{array}\right)+
\left( \begin{array}{c}
\beta\\
2R_2b-R_1\beta
\end{array}\right)Y_t+
\left( \begin{array}{c}
0\\
2\beta
\end{array}\right)r_t\right)dt \\
&\quad+
\left(\begin{array}{c}
\sigma\\
(R_1+2R_2Y_t)\sigma\end{array}\right)dW_t,
\end{align*}
and is therefore clearly a polynomial process with
\[
C_t=\int_0^t \left(\sigma^2\left( \begin{array}{cc}
1 & R_1\\
R_1 & R_1^2\end{array}\right)+\sigma^2\left( \begin{array}{cc}
0 & 2R_2\\
2R_2 & 4R_1R_2\end{array}\right)Y_s+\sigma^2\left( \begin{array}{cc}
0 & 0\\
0 & 4R_2^2\end{array}\right)Y_s^2\right) ds.
\]
\end{example}

\begin{example}[Jacobi process]\label{jacobi}
Another example of a polynomial process is the Jacobi process (see~\cite{gourieroux}), which is the solution of the stochastic differential equation
\begin{align*}
dX_t=-\beta(X_t-\theta)dt+\sigma\sqrt{X_t(1-X_t)}dW_t, \quad X_0=x \in [0,1],
\end{align*}
on $S=[0,1]$, where $\theta \in [0,1]$ and $\beta,\sigma > 0$. 
This example can be extended by adding jumps, where the jump times correspond to those of a Poisson process with intensity $\lambda$ and the jump size is a 
function of the process level. 
Indeed, if a jump occurs, then the process is reflected at $\frac{1}{2}$ so that it remains in the interval $[0,1]$. The extended generator is given by
\[
\mathcal{G}g(x)=-\beta(x-\theta)\frac{dg(x)}{dx}+ \frac{1}{2}\sigma^2(x(1-x))\frac{d^2 g(x)}{dx^2}+\lambda (g(1-x)-g(x)).
\]
In terms of Remark~\ref{rem:example} (ii) (b), we have here $p(x,y)=-2yx+y$ and $\mu(dy)=\lambda\delta_{1}(dy)$.
\end{example}

\begin{example}[Pearson diffusions]
The above example~\ref{jacobi} (without jumps) as well as Ornstein-Uhlenbeck and Cox-Ingersoll-Ross processes, all of them with mean-reverting drift, can be subsumed under the class of so called Pearson diffusions, which are the solutions to SDEs of the form
\[
dX_t=-\beta(X_t-\theta)dt+\sqrt{(a+\alpha_{10}X_t+\alpha_{11}X_t^2)}dW_t, \quad X_0=x,
\]
where $\beta>0 $ and $\alpha_{10},\, \alpha_{11} $ and $a$ are specified such that the square root is well defined. In view of Theorem~\ref{cor:suffcond} it is thus obvious that these processes are polynomial.
Forman and S\o rensen~\cite{forman} give a complete classification of the different types of the Pearson diffusion in terms of their invariant distributions.
\end{example}

\begin{example}[Dunkl process]
The extended generator of the so called Dunkl process (see~\cite{dunkl},~\cite{gallardo}) is given by
\begin{align*}
\mathcal{G}g(x)&=\frac{d^2 g(x)}{dx^2}+\frac{\lambda}{2x^2}\int_{\re} \left(g(x+\xi)-g(x)-\xi\frac{d g(x)}{dx}\right)\delta_{-2x}(d\xi)\\
&=\frac{d^2 g(x)}{dx^2}+\frac{\lambda}{x}\frac{d g(x)}{dx}+\frac{\lambda(g(-x)-g(x))}{2x^2}.
\end{align*}
Since $K(x,d\xi)=\frac{\lambda}{2 x^2}\delta_{-2x}(d\xi)$ and since
\[
\int_{\mathbb{R}} |\xi|^m K(X_t,d\xi)=\frac{\lambda|2X_t|^m}{2X_t^2} =2^{m-1}\lambda |X_t|^{m-2}
\]
for all $m \geq 2$, we derive from Theorem~\ref{cor:suffcond} that the Dunkl process is a polynomial process.
\end{example}

\section{Applications}\label{sec-applications}

By Theorem~\ref{th:characterizationpoly} we know that there exists a linear map $A$ such that moments of $m$-polynomial processes can simply be calculated by computing $e^{tA}$. Indeed, by choosing a basis $\langle e_1,\ldots,e_N\rangle$ of $ \PolmS$ the matrix corresponding to this linear map, which we also denote by $A=(A_{kl})_{k,l=1,\ldots, N}$, can be obtained through
$
\mathcal{A}e_k=\sum_{l=1}^{N} A_{kl}e_l.
$
Writing $f$ as $ \sum_{k=1}^{N} \alpha_k e_k $, we then have
\begin{align}\label{eq:polynomialclaimformula}
P_t f=(\alpha_1,\ldots,\alpha_N)e^{tA}(e_1,\ldots,e_N)',
\end{align}
which means that moments of polynomial processes can be evaluated simply by computing matrix exponentials.

By means of the one-dimensional Cox-Ingersoll-Ross process
\[
dX_t=(b+\beta X_t)dt+\sigma\sqrt{X_t}dW_t, \quad b, \sigma \in \mathbb{R}_+, \quad \beta \in \mathbb{R},
\]
we exemplify how moments of order $m$ can be calculated. Its extended generator is given by
\[
\mathcal{A}g(x)=\frac{1}{2}\sigma^2 x\frac{d^2g(x)}{dx^2}+(b+\beta x)\frac{dg(x)}{dx}.
\]
Applying $\mathcal{A}$ to $(x^0,x^1,\ldots,x^m)$ yields the following $(m+1)\times (m+1)$ matrix
\[
A=\left( \begin{array}{cccccc}
0 & \ldots & & & &  \\
b & \beta & 0 & \ldots & &  \\
0 & 2b+\sigma^2 & 2\beta & 0 & \ldots & \\
0 & 0 & 3b+3\sigma^2 & 3\beta & 0 & \ldots \\
& & & & \ddots & \\
0 & \ldots & & &   mb+\frac{m(m-1)}{2}\sigma^2 & m\beta \\
\end{array}\right).
\]
Hence, $\mathbb{E}_x\left[(X_t)^k\right]=P_t x^k=(0,\ldots,1,\ldots,0)e^{tA}(x^0,\ldots,x^k,\ldots,x^m)'$.

\begin{remark}
\begin{itemize}
\item[(i)]
Note that $A$ is a lower triangular matrix, whose eigenvalues are the diagonal elements. Since in this case they are all distinct, the matrix is diagonalizable. Of course, there are many efficient algorithms to evaluate such matrix exponentials (see for example~\cite{golub},~\cite{moler}).
\item[(ii)]
If $ n > 1 $ one has to apply well-known techniques from linear algebra of polynomials, in order to enumerate efficiently a basis of $ \PolmS $ and to exploit sparsity properties of $ A $, see for instance~\cite{reimer}.
\end{itemize}
\end{remark}

\subsection{Moment estimation - Generalized Method of Moments (GMM)}\label{sec:4.1}

In view of this easy and fast technique of moment calculation for polynomial processes, the Generalized Method of Moments (GMM) can be applied for parameter estimation. If we are given a stationary polynomial process $ X $, then typical functionals applied for parameter estimation are of the form
\[
f(X,\theta)=\left( \begin{array}{c}
X_t^{n_1}X_{t+s}^{m_1}-\mathbb{E}_x[X_t^{n_1}X_{t+s}^{m_1}]\\
\vdots\\
X_t^{n_q}X_{t+s}^{m_q}-\mathbb{E}_x[X_t^{n_q}X_{t+s}^{m_q}]\\
\end{array} \right), \quad  n_i, m_i \in \mathbb{N}, \quad 1\leq i \leq q,
\]
where $\theta$ is the set of parameters to be estimated. This functional is indeed simple to evaluate, since
$
\mathbb{E}_x\left[X_t^{n}X_{t+s}^{m}\right]=\mathbb{E}_x\left[X_t^{n}\,\mathbb{E}_{X_t}\left[X_{s}^{m}\right]\right]
$
can also be computed easily. The usual technology of equating time averages to expectations applies for $ f $ and leads to efficient calibration methods. In the case of one-dimensional jump-diffusions, Zhou~\cite{zhou} already uses this method for GMM estimation.

\subsection{Model calibration}\label{sec:4.1a}

In model calibration -- in contrast to estimation of parameter values of a stationary process from time series data -- parameters are chosen such that derivatives' prices are best explained. Here also the polynomial structure can be very helpful: assume a polynomial process $ X $, where derivatives' prices are known from the market. Derivatives' prices are expectations $\mathbb{E}_x[ f(X_t) ]$ for sufficiently many time points $ t > 0 $ and sufficiently many payoffs $ f $, such that we can estimate the curves $ t \mapsto P_t g(x) $ for today's initial value $ x $ and several polynomials $ g $. In other words we need as many derivatives' prices as necessary to calculate (estimate) the prices of some payoffs, which are polynomials in the underlying.

Having now those curves $ t \mapsto P_t g (x) $ it is often a very easy task to read off the parameter values which explain this curve. This will be worked out in a follow-up paper.

\subsection{Pricing  - Variance reduction}\label{sec:4.2}

The fact that moments of polynomial processes are analytically known also gives rise to new and efficient techniques for pricing and hedging.

Let $X$ be an $m$-polynomial process and $G: S \rightarrow \re^n$ a deterministic bi-measurable map such that the discounted price processes are given through
$
S_t=G(X_t)
$
under a martingale measure. Typically $G = \exp$ if $X$ are log-prices. We denote by $F=\phi(S_T)$ a bounded measurable European claim for some maturity $T > 0$, whose (discounted) price at $t \geq 0$ is given by the risk neutral valuation formula
\[
p^F_t=\mathbb{E}_x\left[\phi(S_T)\big|\mathcal{F}_t\right]=\mathbb{E}_{X_t}\left[(\phi\circ G)(X_T)\right].
\]
Obviously, claims of the form
\begin{align}\label{eq:polyclaims}
F=f\circ G^{-1}(S_T)
\end{align}
for $f \in  \mathcal{P}_m$ are analytically tractable, since we have
\[
p^F_t=\mathbb{E}_x\left[(f \circ G^{-1})(S_T)\big|\mathcal{F}_t\right]=P_{T-t}f(G^{-1}(S_t))=e^{(T-t)A}f(G^{-1}(S_t))
\]
for $0 \leq t\leq T$, where $ A $ is the previously defined linear map on $ \PolmS $. Although claims are in practice not of form~\eqref{eq:polyclaims}, the explicit knowledge of the price of polynomial claims can be used for variance reduction techniques based on control variates. Instead of using the estimator
$
\pi^F_0=\frac{1}{L}\sum_{i=1}^L (\phi\circ G)(X_T^i)
$
in a Monte-Carlo simulation, where $X^1_T,\ldots, X^L_T$ are $L$ samples of $X_T$, we can use
\[
\hat{\pi}^F_0=\frac{1}{L}\sum_{i=1}^L \left((\phi\circ G)(X_T^i)-\left(f(X_T^i)-\mathbb{E}_x\left[f(X_T)\right]\right)\right),
\]
where $f \in \PolmS$ is an approximation of $\phi\circ G$ and serves as control variate. Both estimators are unbiased and the second clearly outperforms the first since
$\Var\left(\hat{\pi}^F_0\right) < \Var\left(\pi^F_0\right)$, where the ratio of the variances depends on the accuracy of the polynomial approximation.

It is worth mentioning that the previous pricing algorithm has also important consequences for hedging, since the Greeks for ``polynomial claims'' $ F = f(X_T) $ can also be calculated explicitly and efficiently: The coefficients of the polynomial $ x \mapsto \mathbb{E}_x[f(X_T)] $ can be computed using matrix exponentiation, taking derivatives of this polynomial is then a simple algebraic operation. To be more precise, the sensitivities of the price process with respect to the factors $X$ can be calculated by
\begin{align*}
\nabla p^F_t=\nabla P_{T-t}f(G^{-1}(S_t))\nabla G^{-1}(S_t).
\end{align*}
Assuming a complete market situation, the claim $ \phi(S_T) = \phi \circ G (X_T) $, can be replicated by a trading strategy $ \eta $, i.e.
\[
\phi(S_T) = \mathbb{E}_x[\phi(S_T)] + \int_0^T \eta_t dS_t.
\]
Similarly the polynomial claim $f(X_t)$ can be replicated by the delta-hedging strategy $\nabla p^F_t$ and we conclude that
\[
\phi(S_T) - f(X_T) = \mathbb{E}_x[\phi(S_T)] - \mathbb{E}_x[f(X_T)] + \int_0^T (\eta - \nabla p^F_t) dS_t.
\]
Therefore, if we assume that $ \phi(S_T) - f(X_T) $ has small variance, then also the stochastic integral representing the difference of the cumulative gains and losses of the two hedging portfolios, namely the one built by the unknown strategy $ \eta $ and the one built by the known strategy $ \nabla p^F_t $, is small.

Concerning the approximation of $\phi \circ G$ by a polynomial, let us consider the case $S_t=G(X_{t,1})$ with $G: \mathbb{R} \rightarrow \mathbb{R}_+$, meaning that we only have one asset which depends on the first component of the polynomial process $X$ as it is usually the case in stochastic volatility models. If the Hamburger moment problem for the law of $X_{T,1}$, say $\mu$, admits a unique solution, then the set of all polynomials is dense in $L^2(\mu)$ (see \cite[Theorem 2.3.3]{akhiezer}) and hence also in $L^1(\mu)$. A sufficient conditions for the uniqueness of a solution to this moment problem is 
\begin{align}\label{eq:momentprob}
\left|\mathbb{E}_x\left[X_{T,1}^k\right]\right|\leq C\frac{k!}{R^k}
\end{align}
for some constants $C >0$ and $R> 0$ (see \cite[Example X.4]{reed}). This condition can be assured by the existence of exponential moments of $X_{T,1}$ around $0$, that is, the moment generating function $\mathbb{E}[e^{uX_{T,1}}]$ is finite for all $u \in (-\varepsilon,\varepsilon)$, which is often satisfied in financial applications.

\begin{figure}
\centering
\includegraphics[width=330pt,height= 200pt]{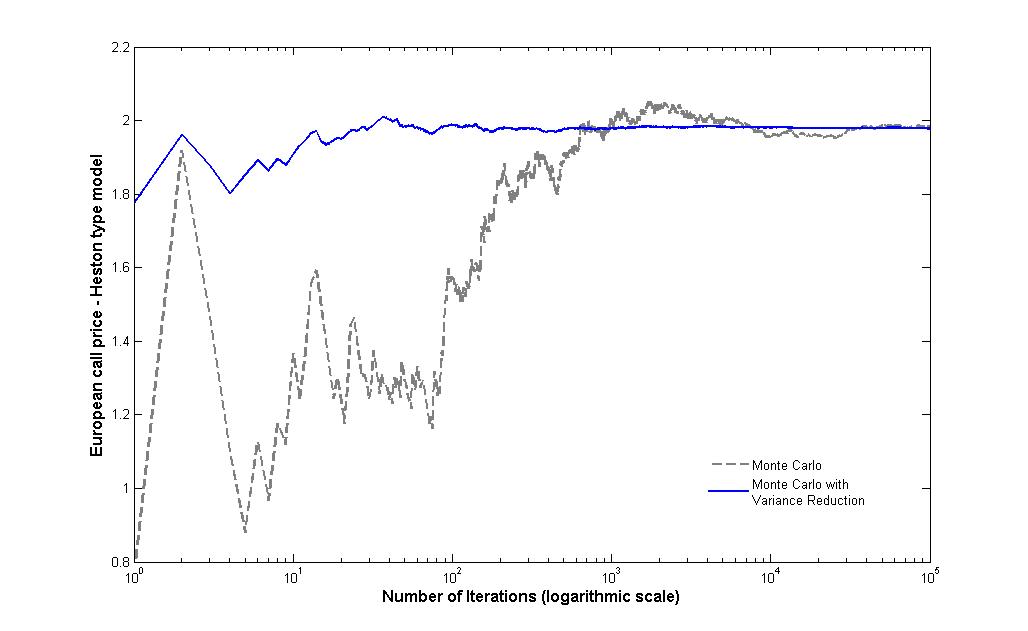}
\caption[Comparison: Monte Carlo simulation, Variance Reduction]{Comparison: Monte Carlo simulation for European call prices with and without variance reduction.}\label{fig:comparison}
\end{figure}

Fur illustratory purposes we have implemented the following affine stochastic volatility model which was initially proposed by Bates~\cite{bates}. The price process is specified as $S_t=S_0e^{X_t}$ with dynamics 
\begin{align*}
d\left(\begin{array}{c}
X_t\\
V_t
\end{array}\right)&=
\left(\begin{array}{c}
r-\frac{V_t}{2}-\lambda V_t\int_{\re}(e^\xi-1)F(d\xi)\\
b-\beta V_t\\
\end{array}\right)dt\\
&\quad+\left( \begin{array}{cc}
\sqrt{V_t} & 0  \\
\sigma\rho\sqrt{V_t} & \sigma\sqrt{1-\rho^2}\sqrt{V_t} \\
\end{array}\right)
\left( \begin{array}{c}
dB_{t,1}\\
dB_{t,2}
\end{array}\right)+\left( \begin{array}{c}
dZ_t\\
0
\end{array}\right),
\end{align*}
where $B$ is a 2-dimensional Brownian motion and $Z$ a pure jump process in $\re$ with jump intensity $\lambda v$ and exponentially distributed jump sizes, i.e.~$F(\xi)=\frac{1}{c}e^{-\frac{\xi}{c}}$, for some parameter $c\in \re_+$. Figure~\ref{fig:comparison} illustrates the comparison between the Monte Carlo simulation for European call prices with and without variance reduction. In this example we use the following parameters:

\begin{table}[h]
\caption{Model parameters}
\label{tab:1}
\begin{tabular}{cccccccccc}
\hline\noalign{\smallskip}
$S_0$ & $V_0$ & Strike & $r$ & $b$ & $\beta$ & $\sigma$ & $\rho$  & $\lambda$ & $c$ \\
\noalign{\smallskip}\hline\noalign{\smallskip}
$10$ & $0.1$ & $9$ & $0.04$ & $0.08$ & $0.7$ & $0.03$ & $0$ & $1.5$ & $0.05$\\
\noalign{\smallskip}\hline
\end{tabular}
\end{table}

The polynomial which we take to approximate the payoff function is of degree $10$ and is chosen such as to minimize the approximation error in a certain interval (depending on the support of the probability distribution).
Concerning computation time, we remark that beside the one-time calculation of the matrix exponential, the only additional computational effort resulting from the use of the control variates is the evaluation of a polynomial in each loop. In our MATLAB code, this causes an increase of computation time of less than $50\%$ per replication. Observing that one can achieve the same accuracy by using $100$ times less replications through the polynomial control variates, the computation time (in our MATLAB implementation) can be decreased by a factor of more than $65$.

Let us finally remark that our variance reduction technique is of particular interest for affine models for which the generalized Riccati ODEs (see~\cite{dfs}) determining the characteristic function cannot be explicitly solved. Moreover, it can also be applied to derivatives involving several assets, provided that their dynamics are described by a polynomial process. This can simply be done by approximating European payoff functions depending on several variables with multivariate polynomials.

\begin{acknowledgements}
All authors gratefully acknowledge the support from the FWF-grant Y 328 (START prize from the Austrian Science Fund) and from ETH foundation. Furthermore the authors are grateful for many comments and valuable suggestions by Chris Rogers, which improved our paper a lot.
\end{acknowledgements}


\end{document}